\documentclass[11pt,a4paper,reqno]{amsart}
\usepackage[T1]{fontenc}
\usepackage{lmodern}
\usepackage{amssymb}
\usepackage{microtype}
\usepackage[top=2.5cm,bottom=2.5cm,left=2cm,right=2cm]{geometry}
\usepackage[hidelinks]{hyperref}

\newtheoremstyle{prelimstyle}
  {6pt plus 2pt minus 2pt}{6pt plus 2pt minus 2pt}
  {\normalfont}{}{\bfseries}{.}{.5em}
  {\thmname{#1}\thmnumber{ #2}\thmnote{ \normalfont#3}}
\theoremstyle{prelimstyle}
\newtheorem{notation}{Notation}[section]
\newtheorem{definition}[notation]{Definition}
\newtheorem{remark}[notation]{Remark}
\newtheorem{example}[notation]{Example}
\theoremstyle{plain}
\newtheorem{theorem}[notation]{Theorem}
\newtheorem{proposition}[notation]{Proposition}
\newtheorem{lemma}[notation]{Lemma}
\newtheorem{corollary}[notation]{Corollary}
\numberwithin{equation}{section}
\DeclareMathOperator{\lcm}{lcm}
\DeclareMathOperator{\Sym}{Sym}
\DeclareMathOperator{\rt}{rt}
\newcommand{\Rees}{\mathcal R}
\makeatletter
\renewcommand\section{\@startsection{section}{1}{\z@}
  {2.5ex \@plus 1ex \@minus .2ex}{1.2ex \@plus .2ex}
  {\normalfont\bfseries}}
\renewcommand\subsection{\@startsection{subsection}{2}{\z@}
  {2ex \@plus 1ex \@minus .2ex}{.8ex \@plus .2ex}
  {\normalfont\bfseries}}
\makeatother

\hypersetup{pdftitle={Linear type and Rees equations of monomial ideals},pdfauthor={Neeraj Kumar}}
\begin{document}
\title[Linear type and Rees equations]{Linear type and Rees equations of monomial ideals}
\author[Neeraj Kumar]{Neeraj Kumar}
\address{Department of Mathematics, Indian Institute of Technology Hyderabad, Kandi, Sangareddy - 502284, INDIA}
\email{neeraj@math.iith.ac.in}

\date{\today}
\begin{abstract}
For a monomial ideal in a polynomial ring over a field, we characterize
linear type by a condition on partitions of its minimal generators. The
same partitions decide whether a new Rees equation is needed in a given
degree, and for ideals with three minimal generators, degree two
already decides linear type. In that case, after the generators are
divided by their greatest common divisor, linear type is equivalent to
an explicit triple of generators being a $d$-sequence. For a
squarefree monomial ideal whose generator graph is a cactus or is
bipartite, linear type is equivalent to the absence of simplicial even
walks in its facet complex. If the generator graph is bipartite and
contains a cycle, the least degree of a minimal nonlinear Rees equation
is half the length of a shortest cycle.
\end{abstract}

\maketitle

\section*{Introduction}

Let $K$ be a field, $S=K[x_1,\ldots,x_n]$, and $I=(f_1,\ldots,f_r)$ a
nonzero proper ideal of $S$. With $t$ an indeterminate over $S$, the
Rees algebra $\Rees(I)=S[It]$ is a quotient of the polynomial ring
$P=S[t_1,\ldots,t_r]$ via the $S$-algebra map
\[
  \varphi:P\longrightarrow\Rees(I)=S[It],
  \qquad t_i\longmapsto f_i t.
\]
The homogeneous elements of $J=\ker\varphi$ are the defining equations
of $\Rees(I)$. We call them \emph{Rees equations} and measure their
degree in $t_1,\ldots,t_r$. The Rees equations of degree one are the
forms $\sum_i a_it_i$ coming from syzygies $\sum_i a_if_i=0$ over $S$.
The ideal $I$ is of \emph{linear type} if these forms generate $J$, or
equivalently, if the canonical map $\Sym_S(I)\to\Rees(I)$ is an
isomorphism.

\vspace{1mm}

How can one recognize linear type directly from the generators? For a
monomial ideal with minimal generators $f_1,\ldots,f_r$, we ask how to
decide from the products of the generators whether a new Rees equation
is needed in a given degree $s\ge2$. The answer is recorded by the
module of effective relations $E_s(I)$, defined in
Remark~\ref{rem:effective}. It vanishes exactly when every Rees equation
of degree $s$ follows from those of smaller degree. For an ideal of a
commutative ring with at least three generators, Planas-Vilanova
characterized $E_s(I)=0$ by two equalities of colon ideals and a
condition on the syzygies among the generators \cite[Theorem~4.7]{PV}.

For monomial ideals the question can be answered by divisibility. Fix
$s\ge2$, split the indices of the minimal generators into two nonempty
sets $A$ and $B$, and let $I_A$ and $I_B$ be the ideals generated by the
corresponding monomials. Suppose that a monomial is divisible by a
product of $s$ generators from $A$ and by a product of $s$ generators
from $B$, repetitions allowed. Must it also be divisible by a product of
$s$ generators that uses both parts? For the given partition, a positive
answer for every such monomial is exactly the containment
\[
                     I_A^s\cap I_B^s\subseteq I_A I_B I^{s-2}.
\]
Theorem~\ref{thm:partition} shows that $E_s(I)=0$ if and only if this
containment holds for every partition. So $I$ is of linear type if and
only if it holds for every partition and every $s\ge2$. The same
partitions also describe the annihilator of $E_s(I)$
(Corollary~\ref{cor:partition-annihilator}).

For an ideal with three minimal monomial generators, every failure of
linear type is already visible in degree two (Theorem~\ref{thm:three}).
A new quadratic equation is needed exactly when the generators can be
ordered so that their exponents strictly increase in one variable and
strictly decrease in another, and Proposition~\ref{prop:three-module}
describes $E_2(I)$ in this case as a direct sum of cyclic $S$-modules.
Minimal Rees equations can nevertheless occur in arbitrarily high
degree: for every $d\ge3$, the ideal $(x^d,x^{d-1}y,y^d)$ needs a
minimal Rees equation of degree $d$ (Example~\ref{ex:three-high}).

\vspace{1mm}

Micali proved that an ideal generated by a regular sequence is of linear
type \cite[Chapter~I, Theorem~1]{Micali}, and Huneke
\cite[Theorem~3.1]{Huneke} and Valla \cite[Theorem~3.15]{Valla}
independently proved the same for ideals generated by $d$-sequences. For
three monomial generators we prove a converse with an explicit choice of
generators (Theorem~\ref{thm:dsequence-three}). Put
$g=\gcd(f_1,f_2,f_3)$ and $m_i=f_i/g$. Then $I$ is of linear type if and
only if
\[
                    m_1-m_3,\qquad m_2-m_3,\qquad m_3
\]
is a $d$-sequence, and in this case the first two terms form a regular
sequence. The division by $g$ cannot be omitted
(Remark~\ref{rem:dsequence-gcd}). When $f_1,f_2,f_3$ have the same
degree, linear type is also equivalent to the vanishing of the
regularity of $\Rees(I)$ in the variables $t_1,t_2,t_3$
(Corollary~\ref{cor:regy-three}).

\vspace{1mm}

For the edge ideal of a connected graph, Villarreal proved that the ideal
is of linear type exactly when the graph is a tree or has a unique
cycle, which is odd \cite[Corollary~3.2]{Villarreal}. Alilooee and Faridi
proved that a squarefree monomial ideal is of linear type if its facet
complex contains no simplicial even walk \cite[Theorem~4.4]{AF}. Their
five-facet example shows why the converse is delicate: four of the
facets form a simplicial even walk, but the quadratic Rees equation of
this walk follows from the linear ones once the fifth facet is added
\cite[Example~4.6]{AF}. In fact, the five-generator ideal is of linear
type by a theorem of Shen \cite[Theorem~4.6]{Shen}. We give the
calculation in Remark~\ref{rem:patch-example}\textup{(i)}.

\vspace{1mm}

We study squarefree monomial ideals through their generator graph
$G(I)$. Its vertices correspond to the minimal generators of $I$, and two
vertices are adjacent when the generators share a variable
(Definition~\ref{def:generator-graph}). For an edge ideal, the vertices
of $G(I)$ are therefore the edges of the original graph. Fouli and Lin
proved linear type when each component of $G(I)$ is a tree or has a
unique cycle, which is odd \cite[Theorem~3.4]{FL}, and Shen described
the Rees equations when $G(I)$ is a single even cycle
\cite[Proposition~4.4]{Shen}.

\vspace{1mm}

We treat two larger classes of generator graphs. If $G(I)$ is a cactus,
that is, each edge lies on at most one cycle, then every even cycle of
length $2s$ gives an indispensable Rees equation of degree $s$, and
these cycle binomials together with the linear Rees equations generate
the Rees ideal (Theorem~\ref{th:cactus}). If $G(I)$ is bipartite and
contains a cycle, the least degree of a minimal nonlinear Rees equation
is half the length of a shortest cycle (Theorem~\ref{th:bip-girth}).
When the shortest cycles have length $2\ell\ge6$, their alternating
binomials are, up to nonzero scalars, exactly the elements of degree
$\ell$ in every minimal binomial generating set of the Rees ideal. In
both classes, linear type is equivalent to the absence of simplicial
even walks in the facet complex (Corollaries~\ref{cor:cactus-walks}
and~\ref{cor:bip-walk}). For complete bipartite generator graphs,
Theorem~\ref{th:bip-complete} counts the minimal Rees equations in each
nonlinear degree. For example, when $G(I)=K_{2,3}$ the least nonlinear
degree is always two, while the number of minimal quadratic equations
can be one, two or three (Example~\ref{ex:bip-two-three}).

\vspace{1mm}

The paper is organized as follows. Section~\ref{sec:preliminaries}
collects the preliminaries. Section~\ref{sec:partitions} proves the
partition criterion. Section~\ref{sec:three} applies it to ideals with
three minimal generators and gives the explicit $d$-sequences, and
Sections~\ref{sec:cactus} and~\ref{sec:bipartite} treat cactus and
bipartite generator graphs. We tested many examples (using AI tools), and every computation stated in
Sections~\ref{sec:three}--\ref{sec:bipartite} was checked over
$\mathbb Q$ with Macaulay2~\cite{GS}.

\section{Preliminaries}\label{sec:preliminaries}

\subsection*{Notation and conventions} Throughout, $K$ is a field and $S=K[x_1,\ldots,x_n]$. We write
$\mathbb N=\{0,1,2,\ldots\}$ and $[m]=\{1,\ldots,m\}$.

\begin{definition}
Let $I=(f_1,\ldots,f_r)$ be a nonzero proper monomial ideal of $S$,
where $f_1,\ldots,f_r$ are its minimal monomial generators, and let $t$
be an indeterminate over $S$. The \emph{Rees algebra} of $I$ is
\[
             \Rees(I)=S[It]=\bigoplus_{s\ge0} I^s t^s\subseteq S[t],
             \qquad I^0=S.
\]
Let $P=S[t_1,\ldots,t_r]$, with one variable $t_i$ for each generator
$f_i$, and consider the surjection of $S$-algebras
\begin{equation}\label{eq:presentation}
        \varphi:P\longrightarrow\Rees(I),\qquad t_i\longmapsto f_i t \quad (1\le i\le r).
\end{equation}
Its kernel $J=\ker\varphi$ is the \emph{Rees ideal} of $I$. We give each
$t_i$ degree one and the elements of $S$ degree zero. Then $P_s$
consists of the forms of degree $s$ in $t_1,\ldots,t_r$, and
$J_s=J\cap P_s$ consists of the Rees equations of degree $s$.

Let $\Sym_S(I)$ be the \emph{symmetric algebra} of the $S$-module $I$.
Since $\varphi(\sum_i b_it_i)=(\sum_i b_if_i)t$, the map
$(b_1,\ldots,b_r)\mapsto\sum_i b_it_i$ identifies the module of
syzygies of $f_1,\ldots,f_r$ with $J_1$. The surjection
$P\to\Sym_S(I)$ sending $t_i$ to $f_i$ in degree one therefore has
kernel $PJ_1$, and $\Sym_S(I)\cong P/PJ_1$ as graded $S$-algebras.
\end{definition}

\begin{definition}
The ideal $I$ is of \emph{linear type} if the canonical map
$\Sym_S(I)\to\Rees(I)$ is an isomorphism. In the presentations above
this map is the quotient map $P/PJ_1\to P/J$, so $I$ is of linear type
if and only if $J=PJ_1$.
\end{definition}

\begin{remark}[{\cite[Lemma~2.1 and Definition~2.2]{PV}}]\label{rem:effective}
For $s\ge2$, the \emph{module of effective relations} of $I$ in degree
$s$ is
\begin{equation}\label{eq:effective}
                        E_s(I)=J_s/P_1J_{s-1}.
\end{equation}
The submodule $P_1J_{s-1}$ consists exactly of the Rees equations of
degree $s$ that follow from those of smaller degree. Indeed, its elements
are sums of linear forms times Rees equations of degree $s-1$, and
conversely $P_{s-k}J_k=P_1(P_{s-k-1}J_k)\subseteq P_1J_{s-1}$ for
$1\le k<s$. So $E_s(I)=0$ exactly when no new Rees equation is needed in
degree $s$.

The \emph{relation type} $\rt(I)$ is the least integer $d\ge1$ such that
$J$ is generated in degrees at most $d$. If $J$ is generated in degrees
at most $d$, then $E_s(I)=0$ for all $s>d$ by the previous paragraph.
Conversely, if these modules vanish, then $J_s=P_1J_{s-1}$ for all
$s>d$, and iterating gives $J_s=P_{s-d}J_d$. In particular, $I$ is of
linear type if and only if $E_s(I)=0$ for every $s\ge2$.
\end{remark}

\subsection*{Rees equations and fibers}

\begin{notation}
Write $\mathbf{t}=(t_1,\ldots,t_r)$. For $u=(u_1,\ldots,u_r)\in\mathbb N^r$,
put $|u|=\sum_i u_i$, $\operatorname{supp}(u)=\{i:u_i>0\}$,
$f^u=\prod_i f_i^{u_i}$ and $\mathbf{t}^u=\prod_i t_i^{u_i}$. For
$A\subseteq[r]$, write
\[
 I_A=(f_i\mid i\in A),\qquad
 f_A=\prod_{i\in A}f_i,\qquad t_A=\prod_{i\in A}t_i .
\]
For $s\ge1$ and a monomial $M\in I^s$, put
\begin{equation}\label{eq:fiber}
 \mathcal F_s(M)=\{u\in\mathbb N^r:|u|=s,\ f^u\mid M\},\qquad
 U_u(M)=\frac{M}{f^u}\mathbf{t}^u\quad (u\in\mathcal F_s(M)).
\end{equation}
\end{notation}

\begin{remark}\label{disc:fiber-grading}
The \emph{fiber} over $Mt^s$ is the set of monomials $a\mathbf{t}^u$ of
$P$ with $\varphi(a\mathbf{t}^u)=Mt^s$. Since
$\varphi(a\mathbf{t}^u)=af^ut^{|u|}$, this means $|u|=s$ and $a=M/f^u$,
so the fiber consists exactly of the terms $U_u(M)$ with
$u\in\mathcal F_s(M)$. Note that different indices give different
terms: if $u\ne v$, then $\mathbf{t}^u\ne\mathbf{t}^v$, even when
$f^u=f^v$.

We say that the monomials in this fiber have \emph{multidegree}
$(M,s)$. The product of monomials of multidegrees $(M,s)$ and $(N,k)$
has multidegree $(MN,s+k)$. Let $P_{(M,s)}$ be the $K$-span of the
fiber and $J_{(M,s)}=J\cap P_{(M,s)}$. Since $\varphi$ maps the
monomials of multidegree $(M,s)$ to $Mt^s$, and distinct monomials of
$S[t]$ are linearly independent over $K$, each multihomogeneous
component of an element of $J$ again lies in $J$. For $s\ge2$, we write
$E_s(I)_M$ for the component of $E_s(I)$ of multidegree $(M,s)$.
\end{remark}

\begin{example}\label{ex:square-fibers}
Let $I=(x,y)^2\subset K[x,y]$, with $f_1=x^2$, $f_2=xy$ and $f_3=y^2$.
For $u=(a,b,c)$ we have $|u|=a+b+c$, $f^u=x^{2a+b}y^{b+2c}$ and
$\mathbf t^u=t_1^at_2^bt_3^c$. For instance, $u=(2,0,0)$ has
$\operatorname{supp}(u)=\{1\}$, and for $A=\{1\}$ we get $I_A=(x^2)$,
$f_A=x^2$ and $t_A=t_1$. We compare two fibers in degree $s=2$.
\begin{enumerate}
\renewcommand{\labelenumi}{\textup{(\roman{enumi})}}
\item Let $M=x^4y$. The only products of two generators that divide $M$
are $f_1^2=x^4$ and $f_1f_2=x^3y$. So
$\mathcal F_2(M)=\{(2,0,0),(1,1,0)\}$, and the two terms of the fiber
are
\[
 U_{(2,0,0)}(M)=\frac{x^4y}{x^4}t_1^2=yt_1^2,\qquad
 U_{(1,1,0)}(M)=\frac{x^4y}{x^3y}t_1t_2=xt_1t_2 .
\]
Their difference is $t_1\ell_{12}$, where $\ell_{12}=yt_1-xt_2$ lies in
$J_1$ because $yf_1=xf_2$. So its class in $E_2(I)=J_2/P_1J_1$ is zero.

\item Let $N=x^2y^2$. Every product of two generators has the same
degree as $N$, so it divides $N$ only if it equals $N$. The
factorizations $N=f_1f_3=f_2^2$ give
$\mathcal F_2(N)=\{(1,0,1),(0,2,0)\}$, and the two terms of the fiber
give $h=t_1t_3-t_2^2\in J_2$. The class of $h$ is not zero. Indeed, let
$\sum_i b_it_i\in J_1$. The degree-two part of
$b_1x^2+b_2xy+b_3y^2=0$ is a relation over $K$ among three distinct
monomials, so each $b_i$ has zero constant term. Hence
$P_1J_1\subseteq(x,y)P$, while $h$ has nonzero image in
$P/(x,y)P=K[t_1,t_2,t_3]$.
\end{enumerate}
\end{example}

\begin{remark}\label{rem:taylor}
The Rees ideal is generated by binomials, as follows.
\begin{enumerate}
\renewcommand{\labelenumi}{\textup{(\roman{enumi})}}
\item Fix $s\ge1$ and a monomial $M\in I^s$. A sum $\sum_u c_uU_u(M)$
with $c_u\in K$ maps to $(\sum_u c_u)Mt^s$, so it is a Rees equation
exactly when $\sum_u c_u=0$. In that case, for any fixed
$v\in\mathcal F_s(M)$,
\[
       \sum_u c_uU_u(M)
       =\sum_{u\ne v}c_u\bigl(U_u(M)-U_v(M)\bigr).
\]
So $J_{(M,s)}$ is spanned over $K$ by the differences of two terms of
the fiber over $Mt^s$.
\item For distinct $u,v\in\mathbb N^r$ with $|u|=|v|=s$, let
$L=\lcm(f^u,f^v)$. The two terms of the \emph{Taylor binomial}
\begin{equation}\label{eq:taylor}
             h_{u,v}=\frac{L}{f^u}\mathbf{t}^u-\frac{L}{f^v}\mathbf{t}^v
\end{equation}
both map to $Lt^s$. If $u,v\in\mathcal F_s(M)$, then $L\mid M$ and
$U_u(M)-U_v(M)=(M/L)h_{u,v}$. Grouping the terms of a Rees equation of
degree $s$ by their images in $S[t]$ writes it as a sum of elements of
$J$ lying in single fibers. By (i), the binomials $h_{u,v}$ therefore
generate $J_s$ as an $S$-module. For $s\ge2$, this description of the
Rees equations is due to Alilooee and Faridi
\cite[Definition~2.2 and Proposition~2.3]{AF}. If all the $f_i$ are
multiplied by the same monomial $w$, then $L$, $f^u$ and $f^v$ are
multiplied by $w^s$, so the binomials~\eqref{eq:taylor} and the Rees
ideal they generate do not change.
\item In degree one, put $g_{ij}=\gcd(f_i,f_j)$. Taking the pairs
$f_i,f_j$ in~\eqref{eq:taylor} shows that $J_1$ is generated over $S$ by
\[
            \frac{f_j}{g_{ij}}t_i-\frac{f_i}{g_{ij}}t_j,
            \qquad 1\le i<j\le r .
\]
\end{enumerate}
\end{remark}

\begin{definition}
A binomial $h\in J$ is \emph{indispensable} if every minimal binomial
generating set of $J$ contains a nonzero scalar multiple of $h$.
\end{definition}

Since distinct monomials of $S[t]$ are linearly independent over $K$,
every nonzero binomial in $J$ is a scalar multiple of a difference of
two monomials in the same fiber. In particular, it is multihomogeneous.

\begin{remark}\label{rem:minimal-fiber}
Write $f_i=x^{\alpha_i}$, and let $e_j$ be the $j$th unit vector of
$\mathbb N^n$. The exponent vectors of the monomials in $\Rees(I)$ form
the semigroup $\Gamma$ generated by $(e_j,0)$ for $1\le j\le n$ and
$(\alpha_i,1)$ for $1\le i\le r$. Thus $\Rees(I)=K[\Gamma]$, and $J$ is
its toric ideal in $P$. Since $\Gamma\subseteq\mathbb N^{n+1}$, we have
$\Gamma\cap(-\Gamma)=\{0\}$, that is, $\Gamma$ is pointed. For
$M=x^\beta$, the multidegree $(M,s)$ of Remark~\ref{disc:fiber-grading}
is the semigroup degree $(\beta,s)$.
\begin{enumerate}
\renewcommand{\labelenumi}{\textup{(\roman{enumi})}}
\item Let $\mathfrak m=(x_1,\ldots,x_n)\subset S$ and
$\mathfrak n=\mathfrak mP+(t_1,\ldots,t_r)P$. Give $x_j$ weight one and
$t_i$ weight $1+\deg f_i$. Then all monomials of multidegree $(M,s)$
have weight $\deg M+s$, and by graded Nakayama's lemma the images of a
minimal multihomogeneous generating set of $J$ form a $K$-basis of
$J/\mathfrak nJ$. Since $(\mathfrak nJ)_s=\mathfrak mJ_s+P_1J_{s-1}$,
the description~\eqref{eq:effective} shows that, for $s\ge2$ and a
monomial $M\in I^s$, the number of elements of multidegree $(M,s)$ in
such a set is
\begin{equation}\label{eq:minimal-fiber}
 \dim_K (J/\mathfrak nJ)_{(M,s)}
   =\dim_K\bigl(E_s(I)/\mathfrak m E_s(I)\bigr)_M .
\end{equation}

\item Charalambous~et~al.\ associate a graph with each
fiber of degree $b\in\Gamma$: two distinct monomials are adjacent if
their difference lies in the ideal generated by the binomials of $J$ of
degrees $a$ with $b-a\in\Gamma\setminus\{0\}$. Choosing binomials that
connect the components of each fiber graph along a spanning tree gives a
minimal binomial generating set, and every such set arises in this way
up to rescaling \cite[Definition~2.1 and Theorems~2.6 and~2.7]{CKT}. For
distinct monomials $U,V$ in one fiber, $U-V$ is indispensable if and
only if the graph of that fiber has exactly the two components $\{U\}$
and $\{V\}$ \cite[Corollary~2.10]{CKT}.

We will often use this in the following form. Suppose that the fiber
over $Mt^s$ consists of two relatively prime monomials $U$ and $V$. By
Remark~\ref{rem:taylor}(i), $J_{(M,s)}=K(U-V)$. No variable of $P$
divides $U-V$, since it would have to divide both terms. So
$(ZJ)_{(M,s)}=0$ for every variable $Z$ of $P$, and hence
$(\mathfrak nJ)_{(M,s)}=0$. Every multiple in this degree of a binomial
of smaller semigroup degree has a nonconstant multiplier and therefore
lies in $\mathfrak nJ$, so it is zero. Hence $U$ and $V$ are isolated,
and $U-V$ is indispensable. For example, in
Example~\ref{ex:square-fibers}(ii) the fiber over $x^2y^2t^2$ consists of
the relatively prime monomials $t_1t_3$ and $t_2^2$, so $t_1t_3-t_2^2$ is
indispensable.
\end{enumerate}
\end{remark}

\section{A partition criterion for linear type}\label{sec:partitions}

Fix a partition of the minimal generators into two parts. To decide
whether a new Rees equation of degree $s$ is needed, we compare products
of $s$ generators taken within one part with products of $s$ generators
that use both parts.

\begin{samepage}
\begin{theorem}\label{thm:partition}
Let $K$ be a field, let $S=K[x_1,\ldots,x_n]$, and let
$f_1,\ldots,f_r$ be the minimal monomial generators of a nonzero proper
monomial ideal $I$ of $S$. For each partition $\{1,\ldots,r\}=A\sqcup B$
into two nonempty sets, let $I_A=(f_i\mid i\in A)$ and
$I_B=(f_j\mid j\in B)$. Fix an integer $s\ge2$. Then $E_s(I)=0$ if and
only if
\begin{equation}\label{eq:partition}
 I_A^s\cap I_B^s\ \subseteq\ I_A I_B I^{s-2}
 \qquad\text{for every such partition}.
\end{equation}
More generally, a monomial $w\in S$ annihilates $E_s(I)$ if and only if
$w(I_A^s\cap I_B^s)\subseteq I_A I_B I^{s-2}$ for every such partition.
Consequently, $I$ is of linear type if and only if~\eqref{eq:partition}
holds in every degree $s\ge2$.
\end{theorem}
\end{samepage}

We call a product of generators \emph{mixed} if it uses at least one
generator from each part. The mixed products of $s$ generators, with
repetitions allowed, generate $I_A I_B I^{s-2}$, and every monomial in
$I_A^s\cap I_B^s$ is divisible by a product of $s$ generators from each
part. So the containment~\eqref{eq:partition} asks whether such a
monomial must also be divisible by a mixed product.

For a fixed $s$, this is a finite test. The monomials $\lcm(m_A,m_B)$,
where $m_A$ is a product of $s$ generators from $A$ and $m_B$ a product
of $s$ generators from $B$, generate $I_A^s\cap I_B^s$: each of them lies
in the intersection, and every monomial of the intersection is divisible
by some $m_A$ and some $m_B$, hence by their least common multiple. It
is therefore enough to find a mixed product dividing each
$\lcm(m_A,m_B)$.

\begin{definition}
For a monomial $M\in I^s$, let $G_s(M)$ be the graph with vertex set
$\mathcal F_s(M)$ from~\eqref{eq:fiber}, in which distinct vertices $u$
and $v$ are adjacent if $\mathbf{t}^u$ and $\mathbf{t}^v$ have a common
variable $t_i$, that is, if
$\operatorname{supp}(u)\cap\operatorname{supp}(v)\ne\varnothing$. Let
$c_s(M)$ be the number of connected components of $G_s(M)$; a graph with
one vertex is connected.
\end{definition}

\begin{samepage}
\begin{lemma}\label{lem:fibercomponents}
For $s\geq2$ and a monomial $M\in I^s$,
\begin{equation}\label{eq:fibercomponents}
 \dim_K E_s(I)_M=c_s(M)-1 .
\end{equation}
More precisely, the Rees equation $U_u(M)-U_v(M)$ follows from those
of smaller degree if and only if $u$ and $v$ lie in the same connected
component of $G_s(M)$.
\end{lemma}
\end{samepage}

\begin{proof}
We first show that the differences $U_u(M)-U_v(M)$ along the edges of
$G_s(M)$ span $(P_1J_{s-1})_{(M,s)}$ over $K$. Suppose that $t_i$
divides both $\mathbf{t}^u$ and $\mathbf{t}^v$, and let $e_i$ be the
$i$th unit vector. For $\nu\in\{u,v\}$ we have $|\nu-e_i|=s-1$ and
$f^{\nu-e_i}=f^\nu/f_i$, which divides $M/f_i$. So
$\nu-e_i\in\mathcal F_{s-1}(M/f_i)$, and dividing $U_\nu(M)$ by $t_i$
gives $U_{\nu-e_i}(M/f_i)$. Hence
\[
 U_u(M)-U_v(M)=
 t_i\bigl(U_{u-e_i}(M/f_i)-U_{v-e_i}(M/f_i)\bigr)\in P_1J_{s-1}.
\]
Conversely, write an element of $(P_1J_{s-1})_{(M,s)}$ as
$\sum_i t_i\Phi_i$ with $\Phi_i\in J_{s-1}$. Since $t_i$ has
multidegree $(f_i,1)$, only the component $H_i$ of $\Phi_i$ of
multidegree $(M/f_i,s-1)$ contributes, and it is zero unless
$f_i\mid M$. By Remark~\ref{rem:taylor}(i), each $H_i$ is a
$K$-linear combination of differences of two terms of the fiber over
$(M/f_i)t^{s-1}$, and multiplying such a difference by $t_i$ gives the
difference along an edge of $G_s(M)$. This proves the claim.

\vspace{1mm}

Now let $C_1,\ldots,C_c$ be the components of $G_s(M)$, and send a
polynomial $\sum_u a_uU_u(M)$ to the vector in $K^c$ whose $j$th entry
is $\sum_{u\in C_j}a_u$. This map kills every edge difference, and its
kernel is exactly the span of the edge differences. Indeed, choose a
vertex $u^{(j)}$ in each $C_j$. If the coefficients over $C_j$ sum to
zero, the part of the polynomial supported on $C_j$ equals
\[
 \sum_{u\in C_j\setminus\{u^{(j)}\}}a_u\bigl(U_u(M)-U_{u^{(j)}}(M)\bigr),
\]
and a path from $u$ to $u^{(j)}$ writes each difference in parentheses
as a sum of edge differences. By Remark~\ref{rem:taylor}(i), the image of
$J_{(M,s)}$ consists of the vectors whose coordinates sum to zero (any
such vector is obtained by placing its entries on the vertices
$u^{(j)}$), and this space has dimension $c-1$. This
proves~\eqref{eq:fibercomponents}. Finally, the images of $U_u(M)$ and
$U_v(M)$ agree exactly when $u$ and $v$ lie in the same component, which
is the last assertion.
\end{proof}

\vspace{1mm}

Diaconis and Sturmfels showed that a set of moves connects every fiber
if and only if the corresponding binomials generate the toric ideal
\cite[Theorem~3.1]{DS}. Lemma~\ref{lem:fibercomponents} is a version of
this for Rees equations of smaller $t$-degree in a fixed fiber. If we
also join $u$ and $v$ whenever $U_u(M)$ and $U_v(M)$ share a variable
$x_j$, then factoring out $x_j$ gives the same argument for
$\mathfrak mJ_s$, and the edge differences span the $(M,s)$-component of
$P_1J_{s-1}+\mathfrak mJ_s=(\mathfrak nJ)_s$. This larger graph
therefore has the same connected components as the graph in
Remark~\ref{rem:minimal-fiber}(ii), and if it has $c$ components, every
minimal multihomogeneous generating set of $J$ has exactly $c-1$
elements of multidegree $(M,s)$ \cite[Theorems~2.6 and~2.7]{CKT}.

For instance, in Example~\ref{ex:square-fibers}(i) the two vertices share
$t_1$ and are adjacent, while in part~(ii) the two vertices over $Nt^2$,
with $N=x^2y^2$, are isolated. Multiplying $N$ by $x$ gives a fiber with
three terms, in which $yt_1t_2$ shares $t_1$ with $xt_1t_3$ and $t_2$
with $xt_2^2$. With $\ell_{23}=yt_2-xt_3$, the resulting path gives
\[
 xh=(xt_1t_3-yt_1t_2)+(yt_1t_2-xt_2^2)
    =-t_1\ell_{23}+t_2\ell_{12},
\]
so $x[h]=0$ in $E_2(I)$. Now regard $I$ as an ideal of $K[x,y,z]$. No
generator involves $z$, so the fiber over $zNt^2$ consists of
$zt_1t_3$ and $zt_2^2$. They share no $t_i$, so $[zh]\ne0$ by
Lemma~\ref{lem:fibercomponents}. The common factor $z$ joins the two
vertices in the larger graph, and indeed $[zh]=z[h]$ lies in
$(x,y,z)E_2(I)$, so this class gives no minimal generator
in~\eqref{eq:minimal-fiber}.

\begin{proof}[Proof of Theorem~\ref{thm:partition}]
Fix $s\ge2$ and a monomial $w\in S$. For every monomial $M\in I^s$ we
have $\mathcal F_s(M)\subseteq\mathcal F_s(wM)$ and $wU_u(M)=U_u(wM)$.
Since the classes of the differences $U_u(M)-U_v(M)$ span $E_s(I)_M$,
Lemma~\ref{lem:fibercomponents} shows that $wE_s(I)\ne0$ if and only if,
for some $M$, two vertices $u,v\in\mathcal F_s(M)$ lie in different
components of $G_s(wM)$. 

\vspace{1mm}

Suppose first that $w(I_A^s\cap I_B^s)\nsubseteq I_A I_B I^{s-2}$ for
some partition. Since the ideals are monomial, we can choose a monomial
$M\in I_A^s\cap I_B^s$ with $wM\notin I_A I_B I^{s-2}$. A product of $s$
generators from each part divides $M$, so $\mathcal F_s(M)$ has a vertex
$u$ supported in $A$ and a vertex $v$ supported in $B$. No vertex of
$G_s(wM)$ is supported on both parts, because its product would be a
mixed product dividing $wM$. Vertices supported in different parts have
disjoint supports, so no edge joins them. Hence $u$ and $v$ lie in
different components of $G_s(wM)$, and $wE_s(I)\ne0$.

\vspace{1mm}

Conversely, suppose that $wE_s(I)\ne0$, and choose $M$ and
$u,v\in\mathcal F_s(M)$ in different components of $G_s(wM)$. Let $C$
be the component containing $u$, and set
\[
 A=\bigcup_{a\in C}\operatorname{supp}(a),\qquad B=[r]\setminus A .
\]
A vertex outside $C$ has support disjoint from $A$, since a common index
would give an edge to a vertex of $C$. So the vertices of $C$ are
supported in $A$ and all other vertices are supported in $B$. In
particular $v$ is supported in $B$, and both parts are nonempty. Since
$f^u$ and $f^v$ divide $M$, we have $M\in I_A^s\cap I_B^s$. On the other
hand, no mixed product of $s$ generators divides $wM$, since it would
give a vertex of $G_s(wM)$ supported on both parts. Hence
$wM\notin I_A I_B I^{s-2}$, and the containment fails.

\vspace{1mm}

If $r=1$, every nonempty fiber has one vertex and there is no partition
to test, so both conditions hold. Taking $w=1$ gives the criterion for
$E_s(I)=0$, and Remark~\ref{rem:effective} gives the statement about
linear type.
\end{proof}

\begin{corollary}\label{cor:partition-annihilator}
Let $I=(f_1,\ldots,f_r)$ be a nonzero proper monomial ideal of
$S=K[x_1,\ldots,x_n]$, where $K$ is a field and the $f_i$ are the
minimal monomial generators of $I$. Put $I_A=(f_i\mid i\in A)$ for
$A\subseteq[r]$. For every $s\ge2$,
\[
 \operatorname{Ann}_S E_s(I)=
 \bigcap_{\substack{[r]=A\sqcup B\\ A,B\ne\varnothing}}
 \bigl(I_A I_B I^{s-2}:(I_A^s\cap I_B^s)\bigr),
\]
where the empty intersection, for $r=1$, is $S$.
\end{corollary}

\begin{proof}
Both sides are monomial ideals. On the right this holds because colons
and intersections of monomial ideals are monomial. On the left, when a
polynomial multiplies a multihomogeneous class of $E_s(I)$, its distinct
terms give products of distinct multidegrees, so the product is zero
exactly when each of these products is zero. It is therefore enough to
test a monomial $w$. Such a $w$ lies in the right-hand side exactly when
$w(I_A^s\cap I_B^s)\subseteq I_A I_B I^{s-2}$ for every partition, and
by Theorem~\ref{thm:partition} this means $wE_s(I)=0$.
\end{proof}

The implication from $E_s(I)=0$ to the containment in
Theorem~\ref{thm:partition} holds for every ideal of a commutative ring.
To prove it, we use Planas-Vilanova's description of effective relations
by first Koszul homology \cite[Corollary~2.7]{PV}.

\begin{proposition}
Let $R$ be a commutative ring, let $I=(f_i\mid i\in\Lambda)$ be an
ideal of $R$, and let $\Lambda=A\sqcup B$ be a partition into two
nonempty sets. Put $I_A=(f_i\mid i\in A)$ and $I_B=(f_j\mid j\in B)$.
For every integer $s\ge2$ there is a surjective $R$-linear map
\[
 E_s(I)\longrightarrow
 \frac{I_A I^{s-1}\cap I_B I^{s-1}}{I_A I_B I^{s-2}}.
\]
If $E_s(I)=0$, then $I_A I^{s-1}\cap I_B I^{s-1}=I_A I_B I^{s-2}$.
In particular, $I_A^s\cap I_B^s\subseteq I_A I_B I^{s-2}$.
\end{proposition}

\begin{proof}
Write $\mathbf{f}t=(f_i t)_{i\in\Lambda}$. By \cite[Corollary~2.7]{PV},
$E_s(I)\cong H_1(\mathbf{f}t;\Rees(I))_s$. After dropping the common
power of $t$, a cycle of degree $s$ is a finitely supported family
$b=(b_i)_{i\in\Lambda}$ with $b_i\in I^{s-1}$ and $\sum_i f_i b_i=0$.
Its partial sum
\[
 \lambda_A(b)=\sum_{i\in A} f_i b_i=-\sum_{j\in B} f_j b_j
\]
lies in $I_A I^{s-1}\cap I_B I^{s-1}$. Conversely, every element $y$ of
this intersection can be written in two ways,
\[
 y=\sum_{i\in A} f_i a_i=\sum_{j\in B} f_j a_j,
 \qquad a_k\in I^{s-1}\quad(k\in\Lambda),
\]
with only finitely many nonzero coefficients, and taking $a_i$ on $A$
and $-a_j$ on $B$ gives a cycle whose partial sum over $A$ is $y$.

Let $e_i$ be the Koszul basis element corresponding to $f_i t$. Since
$e_j\wedge e_i$ has degree two, the boundaries of degree $s$ are
generated, after dropping powers of $t$, by the elements
\[
 c(f_j e_i-f_i e_j),\qquad i\ne j,\quad c\in I^{s-2},
\]
where $I^0=R$. If $i$ and $j$ lie in the same part, the partial sum over
$A$ is zero; if $i\in A$ and $j\in B$, it is $cf_if_j$. These products
generate $I_A I_B I^{s-2}$. So $\lambda_A$ maps the boundaries onto
$I_A I_B I^{s-2}$ and induces the stated surjection on homology. If
$E_s(I)=0$, the target is zero. The last containment follows from
$I_A^s\subseteq I_A I^{s-1}$ and $I_B^s\subseteq I_B I^{s-1}$.
\end{proof}

\begin{corollary}\label{cor:singleton}
Let $I$ be a squarefree monomial ideal of $S$. In
Theorem~\ref{thm:partition} it suffices to consider partitions with at
least two indices in each part. In particular, every squarefree monomial
ideal with at most three minimal generators is of linear type.
\end{corollary}

\begin{proof}
By symmetry we may assume $A=\{i\}$, so that $I_A=(f_i)$. Let $s\ge2$
and $M\in(f_i^s)\cap I_B^s$, and let $f_j$, with $j\in B$, be a factor of
a product of $s$ generators from $B$ that divides $M$. We claim that the
mixed product $f_i^{s-1}f_j$ divides $M$. A variable dividing $f_i$ has
exponent at most $s$ in $f_i^{s-1}f_j$, because the generators are
squarefree, and exponent at least $s$ in $M$, because $f_i^s\mid M$. Any
other variable has the same exponent in $f_i^{s-1}f_j$ as in $f_j$, which
divides $M$. Hence $M\in I_AI_BI^{s-2}$. If $r\le3$, every partition has
a part with one element, so $I$ is of linear type by
Theorem~\ref{thm:partition}.
\end{proof}

The last assertion of Corollary~\ref{cor:singleton} also follows from
work of Alilooee and Faridi \cite[Corollary~3.13 and Theorem~4.4]{AF},
and Fouli and Lin record the three-generator case in
\cite[Remark~3.5]{FL}.

\begin{example}\label{ex:square}
Let $I=(ab,bc,cd,da)\subset K[a,b,c,d]$, with the generators in this
order. Take $s=2$, $A=\{1,3\}$ and $B=\{2,4\}$, so that $I_A=(ab,cd)$ and
$I_B=(bc,da)$. The monomial $M=abcd$ equals both $(ab)(cd)$ and
$(bc)(da)$, so $M\in I_A^2\cap I_B^2$. Each mixed pair of generators
shares a variable, so its product has a square factor and cannot divide
$M$. Hence $M\notin I_A I_B$. In terms of fibers, the fiber over $Mt^2$
consists of the two isolated vertices $t_1t_3$ and $t_2t_4$, so
$t_1t_3-t_2t_4$ has nonzero class in $E_2(I)$. The partitions with a part
of size one give nothing by Corollary~\ref{cor:singleton}, and this
partition detects the quadratic equation.
\end{example}

\section{Ideals with three monomial generators}\label{sec:three}

In this section the ideal $I$ has three minimal monomial generators. In
Section~\ref{subsec:degree-two} we show that degree two decides linear
type, and we describe the module $E_2(I)$. In
Section~\ref{subsec:dsequence} we show that, after the generators are
divided by their greatest common divisor, linear type is equivalent to
an explicit triple of generators being a $d$-sequence.

\subsection{Degree two decides linear type}\label{subsec:degree-two}

For three minimal generators, the partitions in
Theorem~\ref{thm:partition} separate one generator from the other two.
So there are only three partitions to test, and the next theorem shows
that it is enough to test them in degree two.

\begin{samepage}
\begin{theorem}\label{thm:three}
Let $K$ be a field and let $I\subset K[x_1,\ldots,x_n]$ be minimally
generated by three monomials $f_i=\prod_{j=1}^n x_j^{a_{ij}}$.
The following conditions are equivalent.
\begin{enumerate}
\item The ideal $I$ is not of linear type.
\item $I$ has a quadratic Rees equation that does not follow from the
linear ones, that is, $E_2(I)\ne0$.
\item After relabeling the generators, there are distinct indices $p,q$
with
\begin{equation}\label{eq:opposite}
 a_{1p}<a_{2p}<a_{3p},\qquad a_{1q}>a_{2q}>a_{3q}.
\end{equation}
\end{enumerate}
When these inequalities hold, take $A=\{2\}$ and $B=\{1,3\}$, so that
$I_A=(f_2)$ and $I_B=(f_1,f_3)$. Then
$\lcm(f_2^2,f_1f_3)\in(I_A^2\cap I_B^2)\setminus I_A I_B$.
\end{theorem}
\end{samepage}

\begin{proof}
To prove that (1) implies (3), suppose that $I$ is not of linear type.
By Theorem~\ref{thm:partition}, the containment~\eqref{eq:partition}
fails for some partition in some degree $s\ge2$. With three generators,
one part has a single element, so after relabeling $A=\{2\}$ and
$B=\{1,3\}$. Choose a monomial $M\in I_A^s\cap I_B^s$ outside
$I_A I_B I^{s-2}$. There are integers $k,\ell\ge0$ with $k+\ell=s$ such
that
\[
                  f_2^s\mid M,\qquad f_1^k f_3^\ell\mid M .
\]
We first check that $k$ and $\ell$ are both positive. If $\ell=0$, then
$f_1f_2^{s-1}$ divides $\lcm(f_1^s,f_2^s)$, because in each coordinate
one exponent of $f_1$ plus $s-1$ exponents of $f_2$ is at most $s$ times
the larger of the two. Then $f_1f_2^{s-1}$ would be a mixed product
dividing $M$, which is impossible. The same argument with $f_3$ in place
of $f_1$ rules out $k=0$.

\vspace{1mm}

The mixed product $f_2^{s-1}f_3$ does not divide $M$, so there is a
variable $x_p$ whose exponent in $f_2^{s-1}f_3$ is larger than its
exponent in $M$. Since $f_2^s\mid M$, this gives
$(s-1)a_{2p}+a_{3p}>sa_{2p}$, that is, $a_{3p}>a_{2p}$. Comparing with
$f_1^kf_3^\ell\mid M$ and subtracting $a_{3p}$, we also get
\[
 (s-1)a_{2p}>k a_{1p}+(\ell-1)a_{3p}.
\]
If $a_{1p}\ge a_{2p}$, then the right-hand side is at least
$(k+\ell-1)a_{2p}=(s-1)a_{2p}$, because $\ell\ge1$ and $a_{3p}>a_{2p}$.
This contradicts the strict inequality, so $a_{1p}<a_{2p}<a_{3p}$. In
the same way $f_2^{s-1}f_1$ does not divide $M$, and interchanging the
roles of $f_1$ and $f_3$ gives a variable $x_q$ with
$a_{3q}<a_{2q}<a_{1q}$. The two orders show that $p\ne q$, which proves
(3).

\vspace{1mm}

To deduce (2) from (3), label the generators as in~\eqref{eq:opposite},
take $A=\{2\}$ and $B=\{1,3\}$, and put $L=\lcm(f_2^2,f_1f_3)$. Clearly
$L\in I_A^2\cap I_B^2$. At $x_p$, the exponent $a_{2p}+a_{3p}$ of
$f_2f_3$ exceeds both $2a_{2p}$ and $a_{1p}+a_{3p}$, so $f_2f_3$ does not
divide $L$. At $x_q$ the reverse order gives
$a_{1q}+a_{2q}>\max\{2a_{2q},a_{1q}+a_{3q}\}$, so $f_1f_2$ does not
divide $L$ either. Since $f_1f_2$ and $f_2f_3$ generate $I_AI_B$, we get
$L\notin I_A I_B$, and Theorem~\ref{thm:partition} gives $E_2(I)\ne0$.

Finally, (2) implies (1), since $E_2(I)=0$ for every ideal of linear type
by Remark~\ref{rem:effective}.
\end{proof}

\begin{remark}\label{rem:three-classical}
	For an ideal $\mathfrak b$ of a local ring, let $\mu(\mathfrak b)$ be its
	minimal number of generators.
	\begin{enumerate}
		\renewcommand{\labelenumi}{\textup{(\roman{enumi})}}
		\item Fix $p<q$ and put $\mathfrak p=(x_p,x_q)$. In $S_{\mathfrak p}$
		each $f_i$ is a unit multiple of $x_p^{a_{ip}}x_q^{a_{iq}}$, so
		$I_{\mathfrak p}$ is minimally generated by one copy of each of these
		monomials that is minimal under divisibility. Thus
		$\mu(I_{\mathfrak p})=3$ exactly when the pairs $(a_{ip},a_{iq})$ are
		pairwise incomparable in the componentwise order, that is, when the
		first entries are distinct and ordering them increasingly puts the
		second entries in strictly decreasing order. Theorem~\ref{thm:three}
		therefore says that
		\[
		I\text{ is of linear type}
		\quad\Longleftrightarrow\quad
		\mu(I_{(x_p,x_q)})\le2\quad\text{for every }1\le p<q\le n .
		\]
		
		\item Put $g=\gcd(f_1,f_2,f_3)$ and $I'=(f_1/g,f_2/g,f_3/g)$. Dividing by
		$g$ subtracts the same integer from all three exponents in each
		coordinate, so the criterion in (i) is the same for $I$ and $I'$. Since
		the generators of $I'$ have no common variable, $I'$ has height two or
		three, and in height three they form a regular sequence.
		
		Suppose that $\operatorname{ht}I'=2$. At a prime that does not contain
		$I'$, the localization is the unit ideal. A prime of height two
		containing $I'$ is a minimal prime $(x_p,x_q)$, and there $I'$ is a
		complete intersection exactly when it has two minimal generators. At a
		minimal prime of height three, the three generators form a system of
		parameters of a regular local ring, hence a regular sequence. So $I'$ is
		of linear type if and only if it is a complete intersection at each of
		its minimal primes. Simis and Toh\u{a}neanu prove this equivalence for a
		height-two ideal minimally generated by three forms in $K[x,y,z]$, with
		$K$ algebraically closed of characteristic zero
		\cite[Proposition~1.6]{ST}. Here the Taylor resolution gives
		$\operatorname{depth}_{\mathfrak m}(S/I')\ge n-3=\dim(S/I')-1$, where
		$\mathfrak m=(x_1,\ldots,x_n)$, so the implication from complete
		intersection to linear type also follows from
		\cite[Proposition~1.7]{ST}, which Simis and Toh\u{a}neanu attribute to
		Simis and Vasconcelos \cite[Proposition~3.7]{SV}.
	\end{enumerate}
\end{remark}

Theorem~\ref{thm:three} tells us when a new quadratic Rees equation is
needed. To write these equations down, fix $i\in\{1,2,3\}$, let $j$ and
$k$ be the other two indices, and set
\begin{equation}\label{eq:three-Q}
	L_i=\lcm(f_i^2,f_jf_k),\qquad
	C_i=(f_if_j,f_if_k):L_i,\qquad
	Q_i=\frac{L_i}{f_i^2}t_i^2-\frac{L_i}{f_jf_k}t_jt_k .
\end{equation}
Both terms of the Taylor binomial~$Q_i$ map to $L_it^2$, so
$Q_i\in J_2$, and its $x$-multidegree is the exponent vector
$\beta_i\in\mathbb N^n$ of $L_i=x^{\beta_i}$. A monomial $w$ lies in
$C_i$ exactly when $f_if_j$ or $f_if_k$ divides $wL_i$. The next
proposition shows that $C_i$ is the annihilator of $[Q_i]$ in $E_2(I)$
and says when this class is nonzero.

\begin{proposition}\label{prop:three-module}
	Let $K$ be a field and let $I\subset S=K[x_1,\ldots,x_n]$ be minimally
	generated by the three monomials $f_i=\prod_{\nu=1}^n x_\nu^{a_{i\nu}}$,
	$1\le i\le3$. Sending the generator of the $i$th summand to $[Q_i]$
	gives an isomorphism of multigraded $S$-modules
	\begin{equation}\label{eq:three-module}
		\bigoplus_{i=1}^3 (S/C_i)(-\beta_i)\ \simeq\ E_2(I).
	\end{equation}
	For each $i$, let $\{j,k\}=\{1,2,3\}\setminus\{i\}$. The $i$th summand
	is nonzero if and only if there are distinct coordinates
	$p,q\in\{1,\ldots,n\}$ such that
	\[
	a_{jp}<a_{ip}<a_{kp},\qquad a_{jq}>a_{iq}>a_{kq}.
	\]
	Up to nonzero scalars, the $Q_i$ for these indices are exactly the
	quadratic equations in every minimal binomial generating set of the
	Rees ideal $J$.
\end{proposition}

\begin{proof}
By Remark~\ref{rem:taylor}, $J_2$ is generated over $S$ by the Taylor
binomials~\eqref{eq:taylor} of pairs of degree-two monomials in
$t_1,t_2,t_3$. If the two monomials share a variable $t_i$, their
vertices are adjacent in the fiber graph, so the class of the binomial
in $E_2(I)$ vanishes by Lemma~\ref{lem:fibercomponents}. For the pair
$t_i^2,t_j^2$, let $a$ and $b$ be the exponents of a given variable in
$f_i$ and $f_j$. Since $a+b\le\max\{2a,2b\}$, we have
$f_if_j\mid\lcm(f_i^2,f_j^2)$, so the vertex $t_it_j$ joins the two
squares by a path, and again the class vanishes. The only remaining
pairs are $t_i^2,t_jt_k$, which give the $Q_i$. Hence the classes
$[Q_1],[Q_2],[Q_3]$ generate $E_2(I)$.

\vspace{1mm}

To compute the annihilator of $[Q_i]$, fix a monomial $w\in S$ and
consider the class of $wQ_i$. The fiber over $wL_it^2$ contains the
vertices $t_i^2$ and $t_jt_k$. The only possible neighbors of $t_i^2$
are $t_it_j$ and $t_it_k$, and each of them is also adjacent to
$t_jt_k$. Hence $w[Q_i]=0$ exactly when $f_if_j$ or $f_if_k$ divides
$wL_i$, that is, when $w\in C_i$. For a polynomial multiplier it is
enough to check each monomial separately, since distinct monomials
times $Q_i$ have different multidegrees. So $C_i$ is the annihilator of
$[Q_i]$, and since $Q_i$ has $x$-multidegree $\beta_i$, the map
$1+C_i\mapsto[Q_i]$ is a degree-preserving isomorphism from
$(S/C_i)(-\beta_i)$ onto $S\cdot[Q_i]$.

\vspace{1mm}

To see that the sum is direct, fix a multidegree $M$ in the
$x$-variables. If a multiple of $[Q_i]$ is nonzero in degree $M$, then
neither $f_if_j$ nor $f_if_k$ divides $M$. On the other hand, a multiple
of $Q_j$ in degree $M$ exists only if $L_j\mid M$, and $L_j$ is divisible
by $f_if_k$. Similarly, a multiple of $Q_k$ in degree $M$ requires
$f_if_j\mid L_k\mid M$. Hence at most one summand is nonzero in degree
$M$, and on that summand injectivity follows from the annihilator
computation. This proves~\eqref{eq:three-module}.

\vspace{1mm}

Taking $w=1$ above shows that $[Q_i]\ne0$ exactly when neither $f_if_j$
nor $f_if_k$ divides $L_i$. The product $f_if_k$ fails to divide $L_i$
if and only if, in some coordinate $p$, its exponent $a_{ip}+a_{kp}$
exceeds both $2a_{ip}$ and $a_{jp}+a_{kp}$, and these two inequalities
together say $a_{jp}<a_{ip}<a_{kp}$. Interchanging $j$ and $k$, the
product $f_if_j$ fails to divide $L_i$ if and only if
$a_{jq}>a_{iq}>a_{kq}$ for some coordinate $q$. The two strict orders
force $p\ne q$. This proves the criterion.

\vspace{1mm}

For such an $i$, the fiber over $L_it^2$ consists only of the two terms
of $Q_i$. Indeed, the two mixed vertices $t_it_j$ and $t_it_k$ are
absent, and $f_j^2\mid L_i$ would imply $f_if_j\mid L_i$ because
$f_i^2\mid L_i$; the same holds for $f_k^2$. The two terms of $Q_i$
share no variable $t$, and their coefficients are relatively prime by
the definition of the least common multiple, so the two terms are
relatively prime in $P$. By Remark~\ref{rem:minimal-fiber}(ii), $Q_i$ is
indispensable.

\vspace{1mm}

Finally, let $\mathfrak m=(x_1,\ldots,x_n)$ and
$\mathfrak n=\mathfrak mP+(t_1,t_2,t_3)P$. Every binomial in $J$ is
multihomogeneous, so by Remark~\ref{rem:minimal-fiber}(i) the elements of
degree two in any minimal binomial generating set of $J$ map to a
$K$-basis of $(J/\mathfrak nJ)_2$. If $C_i\ne S$, the monomial ideal
$C_i$ is contained in $\mathfrak m$, so $(S/C_i)/\mathfrak m(S/C_i)\cong K$.
Reducing~\eqref{eq:three-module} modulo $\mathfrak m$ gives
\[
 (J/\mathfrak nJ)_2
 =\frac{J_2}{\mathfrak mJ_2+P_1J_1}
 \cong \frac{E_2(I)}{\mathfrak m E_2(I)}
 \cong \bigoplus_{\substack{1\le i\le3\\ C_i\ne S}} K(-\beta_i),
\]
where the image of $Q_i$ spans the $i$th copy of $K$. So every minimal
binomial generating set of $J$ contains exactly
$\#\{i\in\{1,2,3\}\mid C_i\ne S\}$ quadratic equations, and by
indispensability these are scalar multiples of the corresponding $Q_i$.
\end{proof}

\begin{corollary}
For an ideal minimally generated by three monomials in $n$ variables
over a field, the number of quadratic equations in a minimal binomial
generating set of its Rees ideal is at most $\min\{3,\lfloor n/2\rfloor\}$.
For every $n\ge2$ the bound is attained.
\end{corollary}

\begin{proof}
By Proposition~\ref{prop:three-module}, an index $i$ contributes a
quadratic equation exactly when $f_i$ has the strictly middle exponent in
two coordinates, with the other two generators in opposite strict orders
there. A coordinate has at most one strictly middle exponent, so
different contributing indices need disjoint pairs of coordinates.
Hence there are at most $\lfloor n/2\rfloor$ of them, and of course at
most three.

\vspace{1mm}

For sharpness, let $m=\min\{3,\lfloor n/2\rfloor\}$, choose $m$ distinct
indices in $\{1,2,3\}$, and assign a separate pair of variables to each
chosen index $i$. If $j$ and $k$ are the other two indices, let
$(f_j,f_i,f_k)$ have exponents $(0,1,2)$ on the first variable of the
pair and $(2,1,0)$ on the second. Give all generators exponent zero in
the remaining variables. Already on the first pair the three projected
monomials are pairwise incomparable, so the three monomials minimally
generate their ideal. Each chosen index gives a nonzero summand, and the
upper bound is $m$, so there are exactly $m$ minimal quadratic equations.
\end{proof}

\begin{example}\label{ex:three-high}
Let $I=(x^d,x^{d-1}y,y^d)\subset S=K[x,y]$ with $d\ge2$, and take the
generators $f_1,f_2,f_3$ in this order. For $d=2$ this is the ideal of
Example~\ref{ex:square-fibers}. Mui\~nos and Planas-Vilanova give a
minimal set of Rees equations for $(a_1^p,a_2^p,a_1a_2^{p-1})$, where
$a_1,a_2$ is a regular sequence in a Noetherian local ring and $p\ge2$
\cite[Example~5.1]{MPV}. For $a_1=y$, $a_2=x$ and $p=d$ in $S_{(x,y)}$,
their minimal set consists, up to sign, of
\[
 \begin{gathered}
 \ell_{12}=yt_1-xt_2,\qquad \ell_{23}=y^{d-1}t_2-x^{d-1}t_3,\\
 h_s=y^{d-s}t_2^s-x^{d-s}t_1^{s-1}t_3 \qquad(2\le s\le d).
 \end{gathered}
\]
These Rees equations also minimally generate $J$ in $P=S[t_1,t_2,t_3]$.
To see this, put $\mathfrak n=(x,y,t_1,t_2,t_3)P$. The localization map
$J/\mathfrak nJ\to (J/\mathfrak nJ)_{(x,y)}$ is an isomorphism, since
every $c\in S\setminus(x,y)$ acts on $J/\mathfrak nJ$ through its nonzero
constant term. So the classes of the displayed equations form a
$K$-basis of $J/\mathfrak nJ$, and graded Nakayama's lemma applies. For
monomial ideals with three generators in $K[x,y]$, Iglesias~et~al.\ give
an algorithm that computes a minimal
set of Rees equations \cite[Section~3, Theorem~1 and Algorithm~1]{IOSS}.

The fibers show directly why the $h_s$ are needed. Fix $2\le s\le d$ and
put $M_s=\lcm(f_2^s,f_1^{s-1}f_3)=x^{s(d-1)}y^d$. For $|u|=s$, write
$b=u_2+du_3$. Since $f^u=x^{sd-b}y^b$, we have $f^u\mid M_s$ if and only
if $s\le b\le d$. If $u_3=0$, then $b=u_2\ge s$ forces $u_2=s$. If
$u_3\ge1$, the upper bound forces $u_3=1$ and $u_2=0$. Both cases occur,
so
\[
 \mathcal F_s(M_s)=\{(0,s,0),(s-1,0,1)\}.
\]
The fiber over $M_st^s$ thus consists of the two relatively prime terms
of $h_s$, and $h_s$ is indispensable by
Remark~\ref{rem:minimal-fiber}(ii). Since a minimal generating set has
exactly one element in each degree $2\le s\le d$, the $h_s$ are, up to
scalars, all the nonlinear elements of every minimal binomial generating
set of $J$. The last one is $h_d=t_2^d-t_1^{d-1}t_3$. So degree two
detects the failure of linear type, although $\rt(I)=d$ can be
arbitrarily large.

In degree two, only $f_2$ has the middle exponent in both coordinates.
Here $Q_2=h_2$ and $L_2=M_2=x^{2d-2}y^d$ in~\eqref{eq:three-Q}. For a
monomial $w$, we have $f_1f_2\mid wL_2$ if and only if $x\mid w$, and
$f_2f_3\mid wL_2$ if and only if $y\mid w$. So $C_2=(x,y)$, and
Proposition~\ref{prop:three-module} gives
$E_2(I)\cong (S/(x,y))(-\beta_2)$ with $\beta_2=(2d-2,d)$.

Finally, writing the two linear Rees equations in terms of $y$ and $x$
gives
\[
 (\ell_{12},\ell_{23})=(y,x)
 \begin{pmatrix}
 t_1 & y^{d-2}t_2\\
 -t_2 & -x^{d-2}t_3
 \end{pmatrix},
\]
and the determinant of this matrix is $Q_2$. It has bidegree $(d-2,2)$
in the variables $x,y$ and the variables $t_i$, and its class in
$\Sym_S(I)$ is a Sylvester form. Kustin~et~al.\ show that the
component of $\ker(\Sym_S(I)\to\Rees(I))$ of $x,y$-degree $d-2$ is free
of rank one over $K[t_1,t_2,t_3]$, with this class as a basis
\cite[Corollary~2.16 and Remark~2.17]{KPU}.
\end{example}

\subsection{Explicit \texorpdfstring{$d$}{d}-sequences}\label{subsec:dsequence}

Suppose that an ideal $I$ with three minimal monomial generators is of
linear type. We show that, after dividing the generators by their
greatest common divisor, an explicit triple of generators forms a
$d$-sequence. Conversely, such a triple forces linear type.

Recall the definition. Given an ordered minimal generating set
$a_1,\ldots,a_r$ of an ideal, put $\mathfrak{a}_0=(0)$ and
$\mathfrak{a}_i=(a_1,\ldots,a_i)$ for $1\le i\le r$. The sequence
$a_1,\ldots,a_r$ is a \emph{$d$-sequence} if
\[
 \mathfrak{a}_{i-1}:(a_i a_j)=\mathfrak{a}_{i-1}:a_j\qquad(1\le i\le j\le r).
\]
A regular sequence is a $d$-sequence. Huneke \cite[Theorem~3.1]{Huneke}
and Valla \cite[Theorem~3.15]{Valla} proved that an ideal generated by a
$d$-sequence is of linear type; for both facts see also
\cite[Definition~5.5.2 and Corollary~5.5.5]{HS}. For monomial ideals the
monomial generators themselves need not form a $d$-sequence in any order
(Example~\ref{ex:triangle-dsequence}), so we use differences of
monomials instead.

For a monomial $m$ and indices $1\le p<q\le n$, we write $m|_{p,q}$ for
the monomial obtained from $m$ by setting $x_k=1$ for all
$k\notin\{p,q\}$.

\begin{samepage}
\begin{theorem}\label{thm:dsequence-three}
Let $K$ be a field, and let $f_1,f_2,f_3$ be the minimal monomial
generators of an ideal $I$ in $S=K[x_1,\ldots,x_n]$. Put
$g=\gcd(f_1,f_2,f_3)$, $m_i=f_i/g$ and $I'=(m_1,m_2,m_3)$. Then the
following conditions are equivalent.
\begin{enumerate}
\item The ideal $I$ is of linear type.
\item The ordered triple
\[
                  m_1-m_3,\qquad m_2-m_3,\qquad m_3
\]
is a $d$-sequence of generators of $I'$.
\item For all $1\le p<q\le n$, there are indices $i\ne j$ such that
$m_i|_{p,q}$ divides $m_j|_{p,q}$; equality is allowed.
\end{enumerate}
When these conditions hold, $m_1-m_3,\,m_2-m_3$ is a regular sequence,
and the same is true for every labeling of the generators. If the $f_i$
have the same degree, the three generators in (2) are homogeneous.
\end{theorem}
\end{samepage}

Condition (3) is the negation of condition (3) of
Theorem~\ref{thm:three}, written in terms of the $m_i$. By
Remark~\ref{rem:three-classical}(i), it also says that
$\mu(I'_{(x_p,x_q)})\le2$ for all $p<q$.

\begin{proof}
For $\Phi\in P_s$ we have $\Phi(f_1,f_2,f_3)=g^s\Phi(m_1,m_2,m_3)$. Since
$S$ is a domain, $I$ and $I'$ have the same Rees ideal in $P$, so $I$ is
of linear type if and only if $I'$ is.

Put $c=m_3$, $a_1=m_1-c$, $a_2=m_2-c$ and $\mathfrak a=(a_1,a_2)$. Then
$m_1=a_1+c$ and $m_2=a_2+c$, so $a_1,a_2,c$ generate $I'$. They generate
$I'$ minimally, because their classes in $I'/\mathfrak mI'$, where
$\mathfrak m=(x_1,\ldots,x_n)$, are obtained from the basis given by the
classes of $m_1,m_2,m_3$ by an invertible change of basis. Since
$\gcd(m_1,m_2,m_3)=1$, no variable divides all three monomials. The
minimal primes of $I'$ are generated by variables, so
$\operatorname{ht}I'\ge2$.

\emph{(2) implies (1).} This is the theorem of Huneke and Valla.

\vspace{1mm}

\emph{(1) is equivalent to (3).} Fix $p<q$. Dividing by $g$ subtracts
the same integer from the three exponents of $x_p$, and likewise for
$x_q$. So two of the projections $m_i|_{p,q}$ are comparable under
divisibility exactly when the corresponding projections of the $f_i$
are. By Remark~\ref{rem:three-classical}(i), the three projections are
pairwise incomparable if and only if, after relabeling, the opposite
strict orders~\eqref{eq:opposite} hold in the coordinates $p$ and $q$.
Hence (3) fails exactly when condition (3) of Theorem~\ref{thm:three}
holds, and Theorem~\ref{thm:three} shows that (1) and (3) are
equivalent.

\vspace{1mm}

\emph{(1) implies (2).} Assume (1). We first show that the exponent
differences are linearly independent, then that $a_1,a_2$ is a regular
sequence, and finally that $a_1,a_2,c$ satisfies the defining equalities
of a $d$-sequence.

Write $f_i=x^{\alpha_i}$. We show first that $\alpha_1-\alpha_3$ and
$\alpha_2-\alpha_3$ are linearly independent. Every linear Rees
equation $\sum_i b_it_i$ has all
$b_i$ in $\mathfrak m$. Indeed, if $b_i$ had a nonzero constant term,
the monomial $f_i$ would occur in $\sum_i b_if_i=0$ and would have to
cancel against a term of some $b_jf_j$ with $j\ne i$, which would give
$f_j\mid f_i$. So $J=PJ_1\subseteq\mathfrak mP$. Suppose that
$\alpha_1-\alpha_3$ and $\alpha_2-\alpha_3$ were linearly dependent over
$\mathbb Q$. Clearing denominators, we would get integers
$b_1,b_2,b_3$, not all zero, with $b_1+b_2+b_3=0$ and
$\sum_ib_i\alpha_i=0$. Their positive and negative parts give distinct
$u,v\in\mathbb N^3$ with $|u|=|v|$ and $f^u=f^v$. Then
$\mathbf t^u-\mathbf t^v$ lies in $J$ but not in $\mathfrak mP$, a
contradiction. So $\alpha_1-\alpha_3$ and $\alpha_2-\alpha_3$ are
linearly independent. The exponent vectors of the $m_i$ have the same
differences.

\vspace{1mm}

Next we show that $a_1,a_2$ is a regular sequence. Since $m_1\ne m_3$,
the
element $a_1$ is nonzero, hence a nonzerodivisor of the domain $S$. In
the factorial ring $S$, the element $a_2$ is a nonzerodivisor modulo
$(a_1)$ if and only if $a_1$ and $a_2$ have no common prime factor.
Suppose that a prime element $\pi$ divides both. If $\pi$ is a variable,
it divides both terms of $a_1=m_1-m_3$ and both terms of $a_2=m_2-m_3$,
because a variable divides a difference of two distinct monomials only if
it divides both of them. So $\pi$ divides $m_1$, $m_2$ and $m_3$, contrary to
$\gcd(m_1,m_2,m_3)=1$. Otherwise, since $\alpha_1-\alpha_3$ and
$\alpha_2-\alpha_3$ are linearly independent, there is
$\omega\in\mathbb Z^n$ that is orthogonal to $\alpha_1-\alpha_3$ but not
to $\alpha_2-\alpha_3$. Give $x_k$ the weight $\omega_k$. Then $a_1$ is
homogeneous for this $\mathbb Z$-grading, so its factor $\pi$ is
homogeneous too (in a graded domain, the factors of a homogeneous element
are homogeneous, as one sees by comparing highest and lowest components).
The two terms of $a_2$ have different weights. Writing $a_2=\pi h$ and
comparing homogeneous components, we see that $\pi$ divides $m_2$. But a
prime element that divides a monomial is a variable up to a unit, which
is the case already excluded. Hence $a_1,a_2$ have no common prime
factor, and they form a regular sequence, since $\mathfrak a\subseteq
\mathfrak m$.

\vspace{1mm}

Finally we check the defining equalities of a $d$-sequence for
$a_1,a_2,c$. For $i=1$ both sides are zero, since $S$ is a domain. For $i=2$ and
$j\in\{2,3\}$, cancelling the nonzerodivisor $a_2$ modulo $(a_1)$ gives
$(a_1):(a_2a_j)=(a_1):a_j$. It remains to show
$\mathfrak a:c^2=\mathfrak a:c$. Since $a_1,a_2$ is a regular sequence
in the Cohen--Macaulay ring $S$, every associated prime of
$S/\mathfrak a$ is minimal over $\mathfrak a$ and has height two.

\vspace{1mm}

Let $\mathfrak p$ be an associated prime of $S/\mathfrak a$ that contains
$c$. Then $I'=\mathfrak a+(c)\subseteq\mathfrak p$. This is impossible
if $\operatorname{ht}I'=3$. If $\operatorname{ht}I'=2$, then
$\mathfrak p$ is a minimal prime of the monomial ideal $I'$, so
$\mathfrak p=(x_p,x_q)$ for some $p<q$. By (3), choose $i\ne j$ with
$m_i|_{p,q}\mid m_j|_{p,q}$, and write $m_i=\delta m_i'$ and
$m_j=\delta m_j'$, where $\delta=\gcd(m_i,m_j)$. The divisibility means
that neither $x_p$ nor $x_q$ divides $m_i'$. So modulo $\mathfrak p$ the
image of $m_i'$ is a nonzero monomial, while the image of $m_j'$ is zero
or a different monomial, because $m_i\ne m_j$. Hence $m_j'-m_i'$ is a
unit in $S_{\mathfrak p}$. Now $m_j'm_i=m_i'm_j$, and
$m_1\equiv m_2\equiv m_3=c\pmod{\mathfrak a}$, so
$(m_j'-m_i')c\equiv m_j'm_i-m_i'm_j=0\pmod{\mathfrak a}$. Therefore
$c\in\mathfrak aS_{\mathfrak p}$.

\vspace{1mm}

Now let $\eta\in S$ with $\eta c^2\in\mathfrak a$. At an associated
prime $\mathfrak p$ of $S/\mathfrak a$ not containing $c$, the element
$c$ is a unit of $S_{\mathfrak p}$, so $\eta c\in\mathfrak aS_{\mathfrak p}$.
At an associated prime containing $c$, the same holds because
$c\in\mathfrak aS_{\mathfrak p}$. An element of $S/\mathfrak a$ that
vanishes after localizing at every associated prime is zero, so
$\eta c\in\mathfrak a$. This proves $\mathfrak a:c^2\subseteq
\mathfrak a:c$, and the reverse inclusion is clear. So $a_1,a_2,c$ is a
$d$-sequence.

\vspace{1mm}

Conditions (1) and (3) do not depend on the labeling, so the same
argument applies to every labeling of the generators. If the $f_i$ have
the same degree, so do the $m_i$, and their differences are homogeneous.
\end{proof}

When $f_1,f_2,f_3$ have the same degree, Theorem~\ref{thm:dsequence-three}
can be stated in terms of the regularity of the Rees algebra. Suppose
that all $f_i$ have degree $d$. We give $P$ the bigrading with
$\deg x_j=(1,0)$ and $\deg t_i=(0,1)$, and $\Rees(I)$ the bigrading with
$\Rees(I)_{(a,s)}=(I^s)_{a+sd}\,t^s$, so that $\varphi$ is bihomogeneous
of degree $(0,0)$. Following Aramova~et~al.~\cite{ACD} and
R\"omer \cite[Section~1]{Roemer}, we write $\operatorname{reg}_y\Rees(I)$
for the regularity of $\Rees(I)$ in the variables $t_1,t_2,t_3$, that is,
the
largest $b-i$ such that $P(-a,-b)$ occurs in the $i$th module of a
minimal bigraded free resolution of $\Rees(I)$ over $P$.

\begin{corollary}\label{cor:regy-three}
Let $K$ be a field, and let $f_1,f_2,f_3$ be the minimal monomial
generators of an ideal $I$ in $S=K[x_1,\ldots,x_n]$, all of the same
degree. Put $g=\gcd(f_1,f_2,f_3)$ and $m_i=f_i/g$. Then $I$ is of linear
type if and only if $\operatorname{reg}_y\Rees(I)=0$. In this case
$m_1-m_3,\ m_2-m_3,\ m_3$ is a homogeneous generating $d$-sequence of
$(m_1,m_2,m_3)$.
\end{corollary}

\begin{proof}
Put $I'=(m_1,m_2,m_3)$. The $m_i$ have the same degree, and $\Rees(I)$
and $\Rees(I')$ are both isomorphic to $P/J$ as bigraded $P$-modules, so
$\operatorname{reg}_y\Rees(I)=\operatorname{reg}_y\Rees(I')$.

\vspace{1mm}

Suppose that $I$ is of linear type. By Theorem~\ref{thm:dsequence-three},
the three displayed elements form a $d$-sequence, and they are forms of
the same degree. If $K$ is infinite, R\"omer's criterion
\cite[Corollary~3.2(ii)]{Roemer} (see also \cite[Theorem~1.5(2)]{KVSeq})
gives $\operatorname{reg}_y\Rees(I')=0$. If $K$ is finite, we pass to
$K(z)$ for an indeterminate $z$: this flat extension preserves the colon
equalities that define a $d$-sequence, and it does not change the
bigraded Betti numbers of $\Rees(I)$. So $\operatorname{reg}_y\Rees(I)=0$.

\vspace{1mm}

Conversely, the minimal generators of $J$ correspond to the first module
of a minimal bigraded free resolution of $P/J$, so their degrees in
$t_1,t_2,t_3$ are at most $\operatorname{reg}_y\Rees(I)+1$. If
$\operatorname{reg}_y\Rees(I)=0$, then $J$ is generated by elements of
$t$-degree one, that is, $J=PJ_1$, and $I$ is of linear type.
\end{proof}

For an almost complete intersection, such as $I'$ when
$\operatorname{ht}I'=2$, the existence of a generating $d$-sequence was
known before: over an infinite field when the ideal is a complete
intersection at its minimal primes \cite[Proposition~4.10]{JKS}, and in
characteristic zero when the ideal is of linear type and generated in one
degree \cite[Proposition~3.5(2)]{KVACI}. (By
Remark~\ref{rem:three-classical}(ii), $I'$ is a complete intersection at
its minimal primes exactly when $I$ is of linear type.)
Theorem~\ref{thm:dsequence-three} and Corollary~\ref{cor:regy-three} give
the $d$-sequence explicitly, over any field.

The division by $g$ in Theorem~\ref{thm:dsequence-three} cannot be
omitted.

\begin{remark}\label{rem:dsequence-gcd}
In $S=K[x,y,z]$, let
\[
 I=(f_1,f_2,f_3)=(x^2y,xyz,x^2z)=x(xy,yz,xz).
\]
Here $g=x$ and $I'=(xy,yz,xz)$. By Corollary~\ref{cor:singleton}, $I'$
is of linear type, and by Remark~\ref{rem:taylor}(ii), multiplying its
generators by $x$ does not change the Rees ideal, so $I$ is of linear
type as well. Theorem~\ref{thm:dsequence-three} gives the $d$-sequence
$xy-xz,\ yz-xz,\ xz$ of generators of $I'$.

The same differences formed with the original generators are not a
$d$-sequence. Put $a=f_1-f_3=x^2(y-z)$ and $b=f_2-f_3=xz(y-x)$. For
$c\ne0$, the ideal $(a):c$ is generated by $a/\gcd(a,c)$. Since
$\gcd(a,b)=x$ and $\gcd(a,b^2)=x^2$ up to units,
\[
 (a):b=(x(y-z))\subsetneq (y-z)=(a):b^2.
\]
So the ordered sequence $a,b,f_3$ fails the defining colon equality for
$i=j=2$. In fact $a$ and $b$ have the common factor $x$, so they are not
even a regular sequence.
\end{remark}

For squarefree monomials there is a simple test for the monomials
themselves. Let $f_1,f_2,f_3$ be the minimal generators of a squarefree
monomial ideal. By Tang's criterion, they form a $d$-sequence in this
order if and only if $\gcd(f_1,f_2)$ divides $f_3$ \cite[Theorem~2.1]{Tang};
indeed, squarefreeness gives $\gcd(f_i,f_j)=\gcd(f_i,f_j^2)$ for $i<j$,
so only the divisibility condition for the triple remains.

\begin{example}\label{ex:triangle-dsequence}
Huneke and Swanson use the triangle ideal $I=(xy,xz,yz)$ in
$S=K[x,y,z]$ to show that being a $d$-sequence depends on the choice of
generators \cite[Example~5.5.3]{HS}. The ideal is of linear type by
Corollary~\ref{cor:singleton}. In every ordering of the three monomials,
the variable common to the first two is missing from the third, so by
Tang's criterion no ordering of the monomials is a $d$-sequence. For
instance, in the order $xy,xz,yz$,
\[
             (xy):(yz)=(x),\qquad (xy):\bigl((xz)(yz)\bigr)=S.
\]

\begin{samepage}
Here $g=1$, and Theorem~\ref{thm:dsequence-three} gives instead the
ordered generators $a=xy-yz=y(x-z)$, $b=xz-yz=z(x-y)$ and $c=yz$. The
elements $a$ and $b$ are relatively prime, so they form a regular
sequence, and the last colon condition is $(a,b):c^2=(a,b):c$. Put
$\mathfrak{a}=(a,b)$ and $D=(x-z,x-y)$. The ideal $D$ is prime, since
$S/D\cong K[x]$, and $c$ has image $x^2\ne0$ in this quotient. Moreover
$\mathfrak{a}\subseteq D$ and $cD\subseteq\mathfrak{a}$, since
$c(x-z)=za$ and $c(x-y)=yb$. Hence
$D\subseteq\mathfrak{a}:c\subseteq\mathfrak a:c^2\subseteq D:c^2=D$, so
$\mathfrak{a}:c^2=\mathfrak{a}:c=D$, and $a,b,c$ is a $d$-sequence, as
the theorem predicts.\par
\end{samepage}
\end{example}

\begin{remark}\label{rem:dsequence-general}
The proof of Theorem~\ref{thm:dsequence-three} extends to any number
$r\ge2$ of minimal monomial generators, provided that the ideal $I'$,
obtained by dividing the generators by their greatest common divisor,
has height at least $r-1$. For $r=3$ this condition always holds. In
this article we restrict to three generators. For $r\ge4$ the height
condition cannot be dropped, as the next example shows.
\end{remark}

\begin{example}
	In $S=K[x_1,x_2,x_3,x_4]$, let
	\[
	I=(m_1,m_2,m_3,m_4)=(x_2x_3x_4,\ x_1x_3x_4,\ x_1x_2x_4,\ x_1x_2x_3)
	\]
	be the cover ideal of the complete graph $K_4$. Here $r=4$ and
	$\operatorname{ht}I=2<r-1$. Macaulay2~\cite{GS} shows that the Rees
	ideal of $I$ is generated by $x_1t_1-x_4t_4$, $x_2t_2-x_4t_4$ and
	$x_3t_3-x_4t_4$, so $I$ is of linear type. (This can also be checked by
	hand, using $x_im_i=x_1x_2x_3x_4$ for every $i$.) These equations come
	from the linear syzygies $x_im_i=x_4m_4$. Jayanthan and Kumar use the
	linear syzygies of $I^s$ obtained from them to compute the minimal free
	resolutions of the powers of $I$ \cite[proof of Theorem~5]{JKCover}. On
	the other hand, put $a_i=m_i-m_4$ for $i\le3$. The element
	$x_2(x_4-x_1)$ lies in $(a_1):(a_2a_3)$ but not in $(a_1):a_3$, so
	$a_1,a_2,a_3,m_4$ is not a $d$-sequence. Every relabeling of the
	generators comes from a permutation of the variables, so no labeling
	gives a $d$-sequence.
\end{example}

\begin{samepage}
\section{Rees equations for a cactus generator graph}\label{sec:cactus}

The two alternating products of generators around an even cycle give a
Rees equation. When must this equation, up to a nonzero scalar, occur in
every minimal binomial generating set of the Rees ideal?

\begin{definition}\label{def:generator-graph}
Let $f_1,\ldots,f_r$ be the minimal monomial generators of $I$, and put
$F_i=\{x_j:x_j\mid f_i\}$. Following Fouli and Lin
\cite[Construction~2.16]{FL}, the \emph{generator graph} $G(I)$ of $I$
has vertex set $[r]$, and two distinct vertices $i,j$ are adjacent if
and only if $\gcd(f_i,f_j)\ne1$, that is, $F_i\cap F_j\ne\varnothing$.
If the generators $f_i$ are squarefree, their minimality implies that
the $F_i$ are distinct and that none of them is contained in another.
They are therefore the facets of the simplicial complex
$\Delta=\langle F_1,\ldots,F_r\rangle$, which we call the \emph{facet
complex} of $I$.
\end{definition}
\end{samepage}

\begin{remark}\label{disc:generator-graph}
The generator graph appears in the literature under other names. For any
monomial ideal $I$, the correspondence $i\mapsto f_i$ identifies $G(I)$
with the \emph{gcd-graph} of $I$ in the sense of
Berglund~et~al.~\cite[\S\,2.2]{BBH}. Avramov's \emph{GCD graph} has the elements
of the lcm lattice as vertices and joins two distinct monomials when they
are relatively prime \cite[p.~18]{Avramov}, so its induced subgraph on
$f_1,\ldots,f_r$ is the complement of $G(I)$. For squarefree $I$, the
correspondence $i\mapsto F_i$ identifies $G(I)$ with the intersection
graph of the facets of $\Delta$ in the sense of Afkhami and Khosh-Ahang
\cite[p.~77]{AK}; Alilooee and Faridi call it the \emph{line graph}
$L(\Delta)$ \cite[Definition~4.7]{AF}.
\end{remark}

All graphs in this paper are finite and simple. A cycle is a closed walk
of length at least three in which only the first and last vertices
coincide. We call a graph, possibly disconnected, a \emph{cactus} if each
of its edges lies on at most one cycle.

From now on, $I$ is a squarefree monomial ideal of $S$. Some such
assumption is needed: the generator graph of $(x^2,xy,y^2)$ is a path,
yet its Rees ideal requires the quadratic equation $t_1t_3-t_2^2$
(Example~\ref{ex:square-fibers}(ii)).

Let $C=(i_1,\ldots,i_{2k})$ be an even cycle of $G(I)$, and let $A_C$ and
$B_C$ be the sets of indices in its odd and even positions. Define $M_C$
and the alternating binomial $Q_C$ by
\begin{equation}\label{eq:cycle}
 M_C=\lcm(f_{A_C},f_{B_C}),\qquad
 Q_C=\frac{M_C}{f_{A_C}}t_{A_C}-\frac{M_C}{f_{B_C}}t_{B_C}.
\end{equation}
Both terms map to $M_Ct^k$. Changing the starting vertex or the direction
of the cycle can interchange the two terms, so $Q_C$ is determined up to
sign.

\begin{samepage}
\begin{theorem}\label{th:cactus}
Let $I$ be a squarefree monomial ideal of $S=K[x_1,\ldots,x_n]$, where
$K$ is a field. Suppose that each edge of its generator graph $G(I)$ lies
on at most one cycle. Then the Rees ideal of $I$ is
\[
          J=PJ_1+(\,Q_C:C\text{ is an even cycle of }G(I)\,).
\]
Each $Q_C$ is indispensable. For $s\ge2$, the number of minimal
multihomogeneous Rees equations of degree $s$ equals the number of cycles
of length $2s$. In particular, $I$ is of linear type if and only if
$G(I)$ has no even cycle, and its relation type is
\[
    \rt(I)=\max\bigl(\{1\}\cup\{|C|/2:C\text{ is an even cycle of }G(I)\}\bigr),
\]
where $|C|$ denotes the length of $C$.
\end{theorem}
\end{samepage}

The main tool is the following lemma. It also follows from two results
of Alilooee and Faridi. In a simplicial even walk, each facet meets at
least two distinct facets of the opposite list
\cite[Proposition~3.12(i)]{AF}, and if this fails for two lists of
length $s\ge2$, their Taylor binomial lies in
$P_{s-1}J_1+P_1J_{s-1}=P_1J_{s-1}$ \cite[Theorem~4.2]{AF}. (Alilooee and
Faridi divide the generators by their greatest common divisor first;
this does not change the Taylor binomials and can only remove facet
intersections.) If both products divide $M$, the difference of the two
terms in the fiber over $Mt^s$ is a monomial multiple of the Taylor
binomial by Remark~\ref{rem:taylor}(ii). We include the short direct
proof.

\begin{lemma}\label{lem:replace}
Let $I$ be a squarefree monomial ideal of $S$, let $s\ge2$, and let
$M\in I^s$ be a monomial. Take $u,v\in\mathcal F_s(M)$ with disjoint
supports. Color the indices in $\operatorname{supp}(u)$ red and those in
$\operatorname{supp}(v)$ blue, and let $H$ be the bipartite graph on
these indices in which a red and a blue index are adjacent when their
facets meet. If $H$ has a vertex of degree at most one, then
\[
    U_u(M)-U_v(M)\in P_1J_{s-1}.
\]
\end{lemma}

\begin{proof}
By Lemma~\ref{lem:fibercomponents}, it suffices to find a path from $u$
to $v$ in the fiber graph $G_s(M)$. By symmetry, let $j$ be a blue
vertex of degree at most one. Let $i$ be its red neighbor if it has one,
and any red index otherwise. Replace one occurrence of $i$ in $u$ by $j$,
that is, put $u'=u-e_i+e_j$.

\vspace{1mm}

We check that $u'\in\mathcal F_s(M)$. Clearly $|u'|=s$, and
$f^{u'}=(f^u/f_i)f_j$. Only the variables in $F_j\setminus F_i$ gain an
exponent. Each of them is absent from $f^u$, because no red facet other
than $F_i$ meets $F_j$. So it occurs once in $f^{u'}$, by
squarefreeness, and it divides $M$, since $f_j\mid f^v\mid M$. All other
exponents are at most those in $f^u$, and $f^u\mid M$. Hence
$f^{u'}\mid M$.

The monomial $\mathbf{t}^{u'}$ keeps $s-1\ge1$ red factors of
$\mathbf{t}^u$ and has the blue factor $t_j$ of $\mathbf{t}^v$. So $u'$
differs from both $u$ and $v$ and shares an index with each, and
$u,u',v$ is the required path in $G_s(M)$.
\end{proof}

\begin{lemma}\label{lem:cyclefiber}
Let $I$ be a squarefree monomial ideal of $S$ and let $C$ be a chordless
cycle of even length $2k\ge4$ in $G(I)$. Let $W$ be the union of the
facets indexed by $C$. If no facet outside $C$ is contained in $W$, then
the fiber over $M_Ct^k$ consists of the two terms of $Q_C$, and $Q_C$ is
indispensable.
\end{lemma}

\begin{proof}
Two facets of $C$ in positions of the same parity are not consecutive,
so if they met, their edge would be a chord. Thus the facets in $A_C$ are
pairwise disjoint, and so are those in $B_C$. Both alternating products
are therefore squarefree, and together their supports cover $W$, so
$M_C=\prod_{x\in W}x$.

\vspace{1mm}

Take any product of $k$ generators dividing $M_C$. Each factor has its
whole support in $W$, so by hypothesis it comes from $C$. Write $m_j$
for the multiplicity of the generator in position $j$, and set
$m_{2k+1}=m_1$. Positions $j$ and $j+1$ share a variable, whose exponent
in the product is at least $m_j+m_{j+1}$ and whose exponent in $M_C$ is
one. Hence $m_j+m_{j+1}\le1$ for each $j$. On the other hand,
\[
       \sum_{j=1}^{2k}(m_j+m_{j+1})=2\sum_{j=1}^{2k}m_j=2k ,
\]
a sum of $2k$ terms that are each at most one. So every
$m_j+m_{j+1}$ equals one, and the multiplicities alternate between zero
and one. Choosing $m_1=1$ or $m_1=0$ gives exactly the two products
$f_{A_C}$ and $f_{B_C}$, and both divide $M_C$.

The two corresponding terms share no variable $t_i$, and their
coefficients have no common variable, since they are the quotients of a
least common multiple by its two arguments. So the two terms are
relatively prime in $P$, and Remark~\ref{rem:minimal-fiber}(ii) gives
the assertion.
\end{proof}

We use two facts about the blocks of a connected graph: any two distinct
edges of a block lie on a common cycle, and the graph whose vertices are
the blocks and the cut vertices, with incidence as adjacency, is a tree
\cite[p.~61, Lemmas~3.1.3 and~3.1.4]{Diestel}.

\begin{lemma}\label{lem:cactus-end-cycle}
Let $H$ be a finite connected cactus of minimum degree at least two.
Then $H$ contains a cycle $C$ such that at most one vertex of $C$ has a
neighbor outside $C$.
\end{lemma}

\begin{proof}
A block of $H$ that contains a cycle $C$ has no edge outside $C$:
otherwise such an edge and an edge of $C$ would lie on a common second
cycle, contrary to the cactus condition. So each block is a bridge or a
cycle. If $H$ has only one block, the condition on the minimum degree
forces $H$ itself to be a cycle.

Otherwise, choose a leaf of the tree of blocks and cut vertices. Every
cut vertex lies in at least two blocks, so this leaf is an end block,
that is, a block containing exactly one cut vertex of $H$. If it were a
bridge, its other endpoint would have degree one in $H$. Hence the end
block is a cycle, and only its cut vertex can have a neighbor outside it.
\end{proof}

\begin{proof}[Proof of Theorem~\ref{th:cactus}]
Let $J'=PJ_1+(Q_C:C\text{ is an even cycle})$. Each generator of $J'$
lies in $J$. For the reverse inclusion, by Remark~\ref{rem:taylor} it
suffices to show that every difference $U_u(M)-U_v(M)$ lies in $J'$. We
use induction on its degree $s$. The case $s=1$ holds since
$J_1\subseteq J'$. For $s\ge2$, if the two terms have a common factor
$t_i$, we cancel it and apply induction. So we may assume that the
supports of $u$ and $v$ are disjoint, and form the graph $H$ of
Lemma~\ref{lem:replace}. If $H$ has a vertex of degree at most one, the
difference lies in $P_1J_{s-1}\subseteq J'$. 

\vspace{1mm}

Suppose that every vertex of $H$ has degree at least two. Each component
of $H$ is a cactus, since $H$ is a subgraph of $G(I)$, so
Lemma~\ref{lem:cactus-end-cycle} applied to a component of $H$ gives a
cycle $C$ such that at most one vertex of $C$ has a neighbor outside $C$
in $H$. The cycle is even, because red and blue vertices alternate along
it. Interchanging $u$ and $v$ if necessary, we may assume that the
exceptional vertex, if there is one, is red. Then every red neighbor of a
blue vertex of $C$ lies on $C$. 

\vspace{1mm}

Let $A$ be the set of red indices of $C$ and $B$ the set of blue
indices. Let $\nu$ be obtained from $u$ by subtracting one at each index
in $A$, so that $f^u=f^\nu f_A$. We claim that $f^\nu f_B$ divides $M$.
The cycle $C$ has no chord even in $G(I)$, since a chord would create a
second cycle through an edge of $C$, so its blue facets are pairwise
disjoint and every variable of $f_B$ is inserted exactly once. If such a
variable $x$ occurs in $f^u$, then some red facet containing $x$ meets
the blue facet that inserts $x$. That red facet is a neighbor in $H$, so
it lies on $C$ by our choice of colors, and removing $f_A$ removes at
least one copy of $x$. If $x$ does not occur in $f^u$, its new exponent
is one, and $x\mid f^v\mid M$. All other exponents decrease or stay the
same. This proves the claim. 

\vspace{1mm}

Let $u'$ be obtained from $\nu$ by adding one at each index in $B$, so
that $\mathbf{t}^{u'}=\mathbf{t}^\nu t_B$. Both $f^\nu f_A$ and
$f^\nu f_B$ divide $M$, so $M_C=\lcm(f_A,f_B)$ divides $M/f^\nu$. With
the sign of $Q_C$ chosen so that its first term comes from $A$,
\[
       U_u(M)-U_{u'}(M)=\frac{M}{f^\nu M_C}\,\mathbf{t}^\nu Q_C\in J'.
\]
The remaining terms $U_{u'}(M)$ and $U_v(M)$ share the factor $t_j$ for
each $j\in B$. If they are distinct, we cancel one such $t_j$ and apply
induction. This proves $J=J'$. 

\vspace{1mm}

Next we show that $Q_C$ is indispensable for each even cycle $C$. A
facet outside $C$ cannot meet two distinct facets of $C$: otherwise the
two edges through the outside facet, together with either arc of $C$
between their endpoints, would give a second cycle through an edge of
that arc. So no outside facet can lie in the union of the facets of $C$,
since it would then be contained in the single facet of $C$ that it
meets, contrary to minimality. Since $C$ is chordless,
Lemma~\ref{lem:cyclefiber} shows that $Q_C$ is indispensable. 

\vspace{1mm}

Different cycles give different multidegrees $(M_C,|C|/2)$. Indeed, if
two cycles had the same multidegree, the two-term fiber in
Lemma~\ref{lem:cyclefiber} would give the same two monomials in the
$t_i$. Their indices recover the vertex set of the cycle, and a chordless
cycle is the whole induced graph on its vertex set, so the cycles would
coincide. In each cycle multidegree, the class of $Q_C$ in
$J/\mathfrak nJ$ is nonzero by Remark~\ref{rem:minimal-fiber}, and these
classes are linearly independent because their multidegrees are
distinct. By the generation just proved, they span all components of
$J/\mathfrak nJ$ of degree at least two in the $t_i$. So each even cycle
contributes exactly one minimal Rees equation, of degree $|C|/2$. This
proves the count and the statements about relation type and linear
type.
\end{proof}

The computations in Sections~\ref{sec:cactus} and~\ref{sec:bipartite}
are over $K=\mathbb Q$. For each list of generators, $J$ is the kernel
of the map $t_i\mapsto f_it$ in~\eqref{eq:presentation}. To test whether
a Rees equation follows from the linear ones, we reduce it modulo
$PJ_1$, where $J_1$ is obtained from the syzygies of $f_1,\ldots,f_r$.

\begin{example}\label{ex:square-hexagon}
In $S=K[x_1,\ldots,x_{10}]$, let $I=(f_1,\ldots,f_9)$, where
\[
 (f_1,\ldots,f_9)=(x_1x_4x_5x_{10},\,x_1x_2,\,x_2x_3,\,x_3x_4,\,
 x_5x_6,\,x_6x_7,\,x_7x_8,\,x_8x_9,\,x_9x_{10}).
\]
Each variable $x_j$ divides exactly two generators, and it labels the
edge joining their indices in the picture of $G(I)$ below.
\begin{center}
\setlength{\unitlength}{0.85mm}
\begin{picture}(103,47)(0,-2)
\put(6,21){\line(1,1){16}}\put(22,37){\line(1,-1){16}}
\put(38,21){\line(-1,-1){16}}\put(22,5){\line(-1,1){16}}
\put(38,21){\line(1,1){16}}\put(54,37){\line(1,0){27}}
\put(81,37){\line(1,-1){16}}\put(97,21){\line(-1,-1){16}}
\put(81,5){\line(-1,0){27}}\put(54,5){\line(-1,1){16}}
\multiput(22,5)(0,32){2}{\circle*{1.4}}
\multiput(6,21)(32,0){2}{\circle*{1.4}}
\multiput(54,5)(27,0){2}{\circle*{1.4}}
\multiput(54,37)(27,0){2}{\circle*{1.4}}
\put(97,21){\circle*{1.4}}
\put(38,25){\makebox(0,0){$1$}}\put(22,41){\makebox(0,0){$2$}}
\put(2,21){\makebox(0,0){$3$}}\put(22,1){\makebox(0,0){$4$}}
\put(54,41){\makebox(0,0){$5$}}\put(81,41){\makebox(0,0){$6$}}
\put(101,21){\makebox(0,0){$7$}}\put(81,1){\makebox(0,0){$8$}}
\put(54,1){\makebox(0,0){$9$}}
\small
\put(32,32){\makebox(0,0){$x_1$}}\put(12,32){\makebox(0,0){$x_2$}}
\put(12,10){\makebox(0,0){$x_3$}}\put(32,10){\makebox(0,0){$x_4$}}
\put(44,32){\makebox(0,0){$x_5$}}\put(67.5,33.5){\makebox(0,0){$x_6$}}
\put(91,32){\makebox(0,0){$x_7$}}\put(91,10){\makebox(0,0){$x_8$}}
\put(67.5,8.5){\makebox(0,0){$x_9$}}\put(44,10){\makebox(0,0){$x_{10}$}}
\end{picture}
\end{center}
The cycles $C_4=(1,2,3,4)$ and $C_6=(1,5,6,7,8,9)$ meet at vertex $1$.
Their alternating products satisfy $f_1f_3=x_5x_{10}f_2f_4$ and
$f_1f_6f_8=x_1x_4f_5f_7f_9$, and Macaulay2~\cite{GS} gives
$J=PJ_1+(Q_{C_4},Q_{C_6})$, where
\[
 Q_{C_4}=t_1t_3-x_5x_{10}t_2t_4,\qquad
 Q_{C_6}=t_1t_6t_8-x_1x_4t_5t_7t_9.
\]
Here $M_{C_4}=f_1f_3=x_1x_2x_3x_4x_5x_{10}$ and
$M_{C_6}=f_1f_6f_8=x_1x_4x_5x_6x_7x_8x_9x_{10}$. The fiber over
$M_{C_4}t^2$ consists of the two terms of $Q_{C_4}$, the fiber over
$M_{C_6}t^3$ consists of the two terms of $Q_{C_6}$, and in each fiber
the two terms are relatively prime. So the square contributes one
indispensable quadratic equation and the hexagon one indispensable
cubic equation, and $\rt(I)=3$, as Theorem~\ref{th:cactus} predicts. 

\vspace{1mm}

Going around $C_4$ and then $C_6$, starting at vertex $1$, gives a
closed walk of length ten. Its Rees equation of degree five follows from
$Q_{C_4}$ and $Q_{C_6}$, since
\[
 t_1^2t_3t_6t_8-f_1t_2t_4t_5t_7t_9
   =t_1t_6t_8 Q_{C_4}+x_5x_{10}t_2t_4 Q_{C_6},
\]
so its class in $E_5(I)$ is zero.
\end{example}

To relate Theorem~\ref{th:cactus} to the facet complex of $I$, we recall
the definition of a simplicial even walk.

\begin{definition}[{\cite[Definition~3.8]{AF}}]\label{def:simplicial-even-walk}
Let $I$ be a squarefree monomial ideal of $S$ with facet complex
$\Delta=\langle F_1,\ldots,F_r\rangle$. Let $s\ge2$, and consider a
sequence of $2s$ facets, with repetitions allowed,
\[
             F_{i_1},F_{j_1},\ldots,F_{i_s},F_{j_s}.
\]
For each $k\in[r]$, let $u_k$ and $v_k$ be the numbers of occurrences of
$F_k$ in the odd and even positions, respectively. For a vertex $x$ of
$\Delta$, put
\[
          d_u(x)=\sum_{k:\,x\in F_k}u_k,
          \qquad
          d_v(x)=\sum_{k:\,x\in F_k}v_k .
\]
Thus $d_u(x)$ counts the odd positions occupied by a facet containing
$x$, and $d_v(x)$ counts the even positions with this property. Since the
generators are squarefree, these numbers are the exponents of $x$ in
$f^u$ and $f^v$.

The sequence is a \emph{simplicial even walk} if, for every $i$ with
$u_i>0$ and every $j$ with $v_j>0$, both of the following hold.
\begin{enumerate}
\renewcommand{\labelenumi}{\textup{(\roman{enumi})}}
\item There is a vertex $x\in F_i\setminus F_j$ with $d_u(x)\le d_v(x)$.
\item There is a vertex $y\in F_j\setminus F_i$ with $d_v(y)\le d_u(y)$.
\end{enumerate}
In particular, $\operatorname{supp}(u)\cap\operatorname{supp}(v)=\varnothing$:
if $F_k$ occurred in both an odd and an even position, taking $i=j=k$ in
\textup{(i)} would require a vertex in $F_k\setminus F_k$.
\end{definition}

\begin{remark}\label{rem:walk-implications}
For a squarefree monomial ideal $I$ of $S$ with facet complex $\Delta$,
Alilooee and Faridi give the implications
\begin{equation}\label{eq:walk-implications}
 G(I)\text{ has no even cycle}\quad\Longrightarrow\quad
 \Delta\text{ has no simplicial even walk}\quad\Longrightarrow\quad
 I\text{ is of linear type}.
\end{equation}
Under the identification of $G(I)$ with the line graph of $\Delta$ in
Remark~\ref{disc:generator-graph}, the first implication comes from the
proof of \cite[Theorem~4.8]{AF}. It can also be seen directly. Suppose
that $u,v$ give a simplicial even walk. Then each facet $F_i$ with
$i\in\operatorname{supp}(u)$ meets at least two distinct facets indexed
by $\operatorname{supp}(v)$ \cite[Proposition~3.12(i)]{AF}. Indeed, if
$F_i$ met at most one such facet, let $F_j$ be that facet if it exists,
and any facet of the opposite list otherwise. Then every
$x\in F_i\setminus F_j$ has $d_v(x)=0<d_u(x)$, contrary to condition (i)
of Definition~\ref{def:simplicial-even-walk}. Interchanging $u$ and $v$
gives the same statement on the other side; in particular, a walk uses
at least four distinct facets \cite[Corollary~3.13]{AF}. 

\vspace{1mm}

Now form a graph on the two disjoint supports, keeping exactly the edges
of $G(I)$ that join the two sides. Its minimum degree is at least two.
An endpoint of a longest path has a neighbor on the path other than its
predecessor, so the graph contains a cycle. The two sides alternate along
this cycle, so it is an even cycle of $G(I)$. The second implication
in~\eqref{eq:walk-implications} is the even-walk criterion
\cite[Theorem~4.4]{AF}.
\end{remark}

\begin{samepage}
\begin{corollary}\label{cor:cactus-walks}
Let $I$ be a squarefree monomial ideal of $S$ with facet complex
$\Delta$, and suppose that $G(I)$ is a cactus. Then the following
conditions are equivalent.
\begin{enumerate}
\item $I$ is of linear type.
\item $G(I)$ contains no even cycle.
\item The facet complex $\Delta$ contains no simplicial even walk.
\end{enumerate}
\end{corollary}
\end{samepage}  

\begin{proof}
The equivalence of (1) and (2) is part of Theorem~\ref{th:cactus}. By
\eqref{eq:walk-implications}, (2) implies (3), and (3) implies (1).
\end{proof}

\begin{remark}\label{rem:cactus-comparison}
Fouli and Lin proved linear type when each component of the generator
graph is a tree or has a unique cycle, which is odd
\cite[Theorem~3.4]{FL}. Alilooee and Faridi proved the implication from
(2) to (1) for every squarefree monomial ideal \cite[Theorem~4.8]{AF},
but the converse fails without the cactus assumption. In their
Example~4.9, the vertices of the facet complex can be labeled so that
the generators are $x_1x_2x_4$, $x_1x_3x_5$, $x_2x_3x_6$, $x_1x_3x_7$.
The four facets meet pairwise, so the generator graph is $K_4$, which
contains a four-cycle; yet Macaulay2~\cite{GS} gives $J=PJ_1$.

When the whole generator graph is a single even cycle, Shen gives the
presentation $J=PJ_1+(Q_C)$ \cite[Proposition~4.4]{Shen}.
Theorem~\ref{th:cactus} also gives the generators and their degrees when
several cycles occur in the same component.
\end{remark}

\begin{remark}\label{rem:patch-example}
Alilooee and Faridi give an example in which the quadratic Rees equation
of a simplicial even walk follows from the linear ones
\cite[Example~4.6]{AF}. Write $123$ for $\{x_1,x_2,x_3\}$, and similarly
for the other facets, and fix
\[
 F_1=123,\qquad F_2=345,\qquad F_3=567,\qquad F_4=278.
\]
These facets form a simplicial even walk. For the cycle $C=(1,2,3,4)$ of
their generator graph, formula~\eqref{eq:cycle} gives the quadratic
equation $Q_C=x_4x_8t_1t_3-x_1x_6t_2t_4$, in the fiber over $Mt^2$ with
$M=x_1\cdots x_8$.
\begin{enumerate}
\renewcommand{\labelenumi}{\textup{(\roman{enumi})}}
\item Alilooee and Faridi add the facet $F_5=134$. In the generator graph
the square $1,2,3,4,1$ and the triangle $1,2,5,1$ share an edge, so the
cactus hypothesis of Corollary~\ref{cor:cactus-walks} fails. In
$G_2(M)$ the new facet gives the path
\[
 x_4x_8\,t_1t_3 \longleftrightarrow x_2x_8\,t_5t_3
 \longleftrightarrow x_5x_6\,t_5t_4
 \longleftrightarrow x_1x_6\,t_2t_4 ,
\]
whose consecutive terms share $t_3$, $t_5$ and $t_4$, respectively. So
$Q_C\in PJ_1$ by Lemma~\ref{lem:fibercomponents}. In fact
$F_5=(F_1\cup F_2)\setminus(F_3\cup F_4)$ is a patch covering $F_1$ and
$F_2$, so Shen's theorem gives linear type for the five-generator ideal
\cite[Theorem~4.6]{Shen}.

\item Of the $\binom{8}{3}$ three-element subsets of these eight
vertices, four are the given facets, which leaves $52$ choices for a
fifth facet with three vertices. A search in Macaulay2~\cite{GS} over all
$52$ choices, reducing $Q_C$ modulo the linear Rees equations of the
larger ideal, shows that the remainder is zero exactly for $128$, $134$,
$456$ and $678$, and in each of these four cases the full Rees ideal is
$PJ_1$.

\item A similar search in Macaulay2~\cite{GS} over all
$\binom{48}{2}=1128$ unordered
pairs of the remaining facets shows that for exactly $44$ pairs, adding
both facets gives $Q_C\in PJ_1$, but each of these $44$ ideals still
requires a nonlinear Rees equation. For instance, take
$F_5=124$ and $F_6=356$. The pairs of indices
\[
 (1,3)\longleftrightarrow(3,5)\longleftrightarrow(5,6)
       \longleftrightarrow(4,6)\longleftrightarrow(2,4)
\]
give a path in $G_2(M)$, so $Q_C\in PJ_1$. For this pair,
Macaulay2~\cite{GS} gives the Rees ideal
\[
 J=PJ_1+(h),\qquad h=x_8t_1t_3-x_1t_4t_6.
\]
The fiber of $h$ lies over $(M/x_4)t^2$. Since $x_4$ is missing, neither
$F_2$ nor $F_5$ can occur, and among the remaining facets the only
disjoint pairs are $(1,3)$ and $(4,6)$. So the two relatively prime terms
of $h$ form its whole fiber, and $h$ is indispensable.
\end{enumerate}
\end{remark}

Applied in every degree, the following consequence of
Lemma~\ref{lem:replace} recovers the criterion of Alilooee and Faridi
\cite[Theorem~4.8]{AF}.

\begin{corollary}\label{cor:even-cycle-degree}
Let $I$ be a squarefree monomial ideal of $S$ and let $s\ge2$. If $G(I)$
has no even cycle of length at most $2s$, then $E_s(I)=0$.
\end{corollary}

\begin{proof}
Suppose that $E_s(I)\ne0$. By Remark~\ref{rem:taylor}(i), some
difference $U_u(M)-U_v(M)$ of two terms of a fiber has nonzero class in
$E_s(I)$. The two terms share no $t_i$, since factoring out $t_i$ would
leave a Rees equation of degree $s-1$.

On the disjoint supports of $u$ and $v$, form the graph $H$ of
Lemma~\ref{lem:replace}. By that lemma, every vertex of $H$ has at least
two neighbors, so $H$ contains a cycle. This cycle is even, since each
edge of $H$ joins the two supports, and its length is at most
$|\operatorname{supp}(u)|+|\operatorname{supp}(v)|\le2s$. It is also a
cycle of $G(I)$.
\end{proof}

\section{Degrees of Rees equations for bipartite generator graphs}\label{sec:bipartite}

Let $I$ be a squarefree monomial ideal of $S$ whose generator graph
$G(I)$ is bipartite. If $G(I)$ contains a cycle, what is the least degree
of a minimal nonlinear Rees equation? We show that it is half the length
of a shortest cycle.

Fix a bipartition $V_1\sqcup V_2$ of the vertices of $G(I)$; one part may
be empty when the graph has no edges. Generators in the same part have
disjoint supports, so a product of distinct generators from one part is
squarefree. For nonempty subsets $A\subseteq V_1$ and $B\subseteq V_2$ of
the same size, put
\begin{equation}\label{eq:bip-binomial}
 M_{A,B}=\lcm(f_A,f_B),\qquad
 Q_{A,B}=\frac{M_{A,B}}{f_A}t_A-\frac{M_{A,B}}{f_B}t_B .
\end{equation}
Here each selected generator is used once.

\begin{samepage}
\begin{theorem}\label{th:bip-upper}
Let $K$ be a field and let $I\subset S=K[x_1,\ldots,x_n]$ be a
squarefree monomial ideal. Suppose that $V_1\sqcup V_2$ is a bipartition
of $G(I)$. Then the Rees ideal of $I$ is
\[
                         J=PJ_1+\bigl(Q_{A,B}\bigr),
\]
where $A\subseteq V_1$ and $B\subseteq V_2$ range over the pairs with
$|A|=|B|\ge2$ for which every vertex of $A\cup B$ has at least two
neighbors in $A\cup B$. In particular, with $m=\min\{|V_1|,|V_2|\}$, we
have $\rt(I)\le\max\{1,m\}$.
\end{theorem}
\end{samepage}

\begin{proof}
Let $J'$ be the ideal on the right. Every displayed generator is a Rees
equation, and we show by induction on $s$ that $J_s\subseteq J'$. This
holds in degree one. For $s\ge2$, fix a monomial $M\in I^s$ and two
terms $U_u(M),U_v(M)$ of its fiber; by Remark~\ref{rem:taylor} it is
enough to treat their difference. If the two terms have a common factor
$t_i$, the difference lies in $P_1J_{s-1}$, and induction applies. So we
may assume that $u$ and $v$ have disjoint supports. Let $H$ be the graph
of Lemma~\ref{lem:replace} on these supports. If $H$ has a vertex of
degree at most one, the difference lies in $P_1J_{s-1}$ by that lemma,
and again induction applies. In particular, every omitted $Q_{A,B}$ of
degree $s$ lies in $P_1J_{s-1}$, since for its two terms $H$ is the graph
induced on $A\cup B$.

We may now assume that every vertex of $H$ has at least two neighbors.
An edge of $H$ joins an index of $u$ and an index of $v$ lying in
opposite parts of $G(I)$. After interchanging the parts if necessary,
the sets $X=\operatorname{supp}(u)\cap V_1$ and
$Y=\operatorname{supp}(v)\cap V_2$ are both nonempty, and by
interchanging the two terms and the two parts together we may assume
$|X|\le |Y|$. Choose a set $B$ of $|X|$ distinct indices in $Y$. In the
product $f^u$, remove one copy of each generator indexed by $X$ and
insert one copy of each generator indexed by $B$, and let $u'$ record the
multiplicities in the resulting product of $s$ factors. 

\vspace{1mm}

We check that $f^{u'}$ divides $M$. The generators indexed by $B$ have
disjoint supports, so the insertions add at most one copy of each
variable. Let $x$ divide an inserted generator $f_j$. Then $j$ is not in
$\operatorname{supp}(u)$, and no other generator indexed by $V_2$
contains $x$. So if $x$ occurs in $f^u$, it occurs in a generator
indexed by $X$, and removing one copy of that generator offsets the copy
of $x$ that was inserted. If $x$ does not occur in $f^u$, its new
exponent is one and $x\mid f_j\mid f^v\mid M$. All other exponents stay
the same or decrease. So $u'\in\mathcal F_s(M)$.

\vspace{1mm}

If $|X|<s$, the new product keeps a factor of $u$ and has a factor of
$v$, so $u'$ joins $u$ to $v$ in $G_s(M)$, and the difference lies in
$P_1J_{s-1}$ by Lemma~\ref{lem:fibercomponents}. If $|X|=s$, then $u$
uses exactly the $s$ distinct generators indexed by $X$, and $v$ uses
exactly the $s$ distinct generators indexed by $B$. In this case $H$ is
the graph induced by $G(I)$ on $X\cup B$, since every edge on these
vertices joins $X$ to $B$. Every vertex has at least two neighbors in
this graph, so $Q_{X,B}$ is one of the generators of $J'$, and the
difference is $(M/M_{X,B})Q_{X,B}$ up to sign. This completes the
induction. Finally, each $Q_{A,B}$ in the formula has degree
$|A|=|B|\le m$.
\end{proof}

If $G(I)$ is disconnected, the bipartition is not unique. The least value
of $m$ over all bipartitions is the sum, over the connected components,
of the size of the smaller side: each of $V_1$ and $V_2$ contains at
least that many vertices, and putting the smaller side of each component
into $V_1$ attains this value. An isolated vertex contributes zero.

\vspace{1mm}

For the square $I=(ab,bc,cd,da)$ of Example~\ref{ex:square}, the two
parts are $\{1,3\}$ and $\{2,4\}$, and they are the only possible choice
of $A,B$ of size at least two in Theorem~\ref{th:bip-upper}.
Macaulay2~\cite{GS} gives $J=PJ_1+(t_1t_3-t_2t_4)$, so the quadratic
equation detected by the partition in Example~\ref{ex:square} is the only
additional generator needed.

\begin{samepage}
\begin{theorem}\label{th:bip-girth}
Let $I$ be a squarefree monomial ideal of $S=K[x_1,\ldots,x_n]$, where
$K$ is a field, and suppose that its generator graph $G(I)$ is
bipartite. If $G(I)$ is a forest, then $I$ is of linear type. If the
shortest cycle of $G(I)$ has length $2\ell$, then
\begin{equation}\label{eq:bip-first}
 \min\{s\ge2:E_s(I)\ne0\}=\ell .
\end{equation}
If moreover $\ell\ge3$, then the minimal binomial Rees equations of
degree $\ell$ are, up to nonzero scalars, the alternating binomials $Q_C$
of the cycles $C$ of length $2\ell$, and each of them is indispensable.
\end{theorem}
\end{samepage}

\begin{proof}
If $G(I)$ is a forest, Corollary~\ref{cor:even-cycle-degree} gives
$E_s(I)=0$ for every $s\ge2$, so $I$ is of linear type. If the shortest
cycle has length $2\ell$, the same corollary gives $E_s(I)=0$ for
$2\le s<\ell$. It remains to show that $E_\ell(I)\ne0$, together with
the description of the minimal equations when $\ell\ge3$. 

\vspace{1mm}

Suppose first that $\ell\ge3$, and let $C$ be a cycle of length $2\ell$.
A chord of $C$ would give a shorter cycle, so $C$ is chordless. Suppose
that a facet outside $C$ were contained in the union of the facets of
$C$. Since no minimal facet is contained in another, it would meet at
least two facets of $C$. The shorter arc of $C$ between two such facets
has length at most $\ell$, and adding the two edges through the outside
facet gives a cycle of length at most $\ell+2<2\ell$, contrary to the
choice of $C$. So Lemma~\ref{lem:cyclefiber} applies: the fiber over
$M_Ct^\ell$ consists of the two relatively prime terms of $Q_C$, and
$Q_C$ is indispensable. In particular $E_\ell(I)\ne0$. 

\vspace{1mm}

To see that no other binomials are needed in degree $\ell$, take a fiber
difference with nonzero class in $E_\ell(I)$. Its supports are disjoint,
and by Lemma~\ref{lem:replace} its graph $H$ has minimum degree at least
two, so $H$ contains a cycle. Since $H$ has at most $2\ell$ vertices and
every cycle has at least $2\ell$ vertices, this cycle uses every vertex
of $H$. So each term uses $\ell$ distinct indices once, and the indices
of the two terms alternate around the cycle. Hence the difference is
$(M/M_C)Q_C$ up to sign, and its class in $J/\mathfrak nJ$ is zero if
$M_C$ properly divides $M$. So the binomials that can occur in a minimal
generating set in degree $\ell$ are the $Q_C$ themselves, and each of
them occurs by indispensability. Different shortest cycles give
different multidegrees, since the two-term fiber recovers the vertex set
of the cycle and this set induces the chordless cycle. 

\vspace{1mm}

Finally, we prove $E_2(I)\ne0$ when $\ell=2$. Among all four-cycles,
choose one for which the union $W$ of its facets is minimal under
inclusion. Write $A_1,A_2$ for its facets in one part and $B_1,B_2$ for
its facets in the other. Both alternating pairs divide
$M=\prod_{x\in W}x$. If no outside facet is contained in $W$, the
four-cycle has no chord, and Lemma~\ref{lem:cyclefiber} gives an
indispensable quadratic equation, whose class in $E_2(I)_M$ is nonzero by
Remark~\ref{rem:minimal-fiber}. 

\vspace{1mm}

Otherwise, after interchanging the parts, there is an outside facet
$F\subseteq W$ in the part of $A_1$ and $A_2$. It is disjoint from both,
so it lies in $B_1\cup B_2$, and it meets both $B_1$ and $B_2$, since it
is contained in neither. For $i=1,2$, the facets $F,B_1,A_i,B_2$ form a
four-cycle whose union lies in $W$, so by minimality both unions equal
$W$. A variable in $W\setminus(B_1\cup B_2)$ would lie outside $F$, and
hence in both $A_1$ and $A_2$, which are disjoint. Consequently
$W=B_1\sqcup B_2$. 

\vspace{1mm}

Let $b_1,b_2$ be the monomials with supports $B_1,B_2$, so that
$M=b_1b_2$. A pair of generators dividing $M$ and containing $b_1$ must
use a generator dividing $b_2$, which is $b_2$ itself by minimality.
Likewise, a pair containing $b_2$ is $b_1b_2$. So the vertex for $b_1b_2$
is isolated in $G_2(M)$. The pair from $A_1,A_2$ gives another vertex, so
Lemma~\ref{lem:fibercomponents} gives $E_2(I)_M\ne0$. This
proves~\eqref{eq:bip-first} also for $\ell=2$.
\end{proof}

\begin{example}\label{ex:bip-hexagon}
Keep only the six generators on the hexagon $C_6$ of
Example~\ref{ex:square-hexagon}, with their original indices. They
generate the ideal
\[
 I_{C_6}=(f_1,f_5,f_6,f_7,f_8,f_9)
 =(x_1x_4x_5x_{10},x_5x_6,x_6x_7,x_7x_8,x_8x_9,x_9x_{10}),
\]
whose generator graph is $C_6=(1,5,6,7,8,9)$. With
$P=S[t_1,t_5,\ldots,t_9]$, Macaulay2~\cite{GS} gives the Rees ideal
\[
 J=PJ_1+(Q_{C_6}),\qquad Q_{C_6}=t_1t_6t_8-x_1x_4t_5t_7t_9,
\]
so there is no minimal quadratic equation. To see why the cubic equation
$Q_{C_6}$ is needed, put $M=M_{C_6}=x_1x_4x_5x_6x_7x_8x_9x_{10}$. A
product of three generators dividing this squarefree monomial must use
three distinct, pairwise nonadjacent vertices of the hexagon. Only the
two alternating triples are possible, and their products are
$f_1f_6f_8=M$ and $f_5f_7f_9=M/(x_1x_4)$. So the fiber over $Mt^3$
consists of the two relatively prime terms of $Q_{C_6}$, and $Q_{C_6}$ is
indispensable. The least degree of a minimal nonlinear Rees equation is
three, half the length of the hexagon, as in Theorem~\ref{th:bip-girth}.
Shen gives this presentation, by the linear Rees equations and the
alternating binomial, for any single even cycle
\cite[Proposition~4.4]{Shen}.
\end{example}

\begin{corollary}\label{cor:bip-walk}
Let $I$ be a squarefree monomial ideal of $S$ with facet complex
$\Delta$, and suppose that its generator graph $G(I)$ is bipartite. Then
the following conditions are equivalent.
\begin{enumerate}
\item $I$ is of linear type.
\item $G(I)$ is a forest.
\item $\Delta$ contains no simplicial even walk.
\end{enumerate}
\end{corollary}

\begin{proof}
Suppose that $G(I)$ contains a cycle. Since $G(I)$ is bipartite, a
shortest cycle has even length $2\ell$ with $\ell\ge2$, and
Theorem~\ref{th:bip-girth} gives $E_\ell(I)\ne0$, so $I$ is not of
linear type. This proves $(1)\Rightarrow(2)$. A forest has no even
cycle, so the two implications in Remark~\ref{rem:walk-implications}
give $(2)\Rightarrow(3)\Rightarrow(1)$.
\end{proof}

Deleting $f_1$ from Example~\ref{ex:bip-hexagon} leaves the ideal
$(f_5,\ldots,f_9)$, whose generator graph is a path, and for this ideal
Macaulay2~\cite{GS} gives $J=PJ_1$. Adding $f_1$ back closes the path
into the hexagon and introduces the indispensable cubic equation
$Q_{C_6}$, which illustrates Corollary~\ref{cor:bip-walk}. For a
generator graph that is not bipartite, see
Remark~\ref{rem:patch-example}: there the graph contains a triangle, and
the quadratic Rees equation of a simplicial even walk follows from the
linear ones. 

\vspace{1mm}

We now determine the minimal Rees equations in all nonlinear degrees
when $G(I)=K_{p,q}$ with $2\le p\le q$. We call the generators in the
part of size $p$ the rows $a_1,\ldots,a_p$, and those in the other part
the columns $b_1,\ldots,b_q$. Put
\begin{equation}\label{eq:bip-column-factors}
 a=\prod_{i=1}^p a_i,\qquad
 b'_j=\frac{b_j}{\gcd(b_j,a)}\quad(1\le j\le q).
\end{equation}
The variables dividing $b'_j$ are exactly those that divide $b_j$ and no
other generator. So each $b'_j$ is relatively prime to $a$, and
$b'_1,\ldots,b'_q$ are pairwise relatively prime. For row indices
$A\subseteq[p]$ and column indices $B\subseteq[q]$ with $|A|=|B|\ge1$,
we write $M_{A,B}$ and $Q_{A,B}$ for the monomial and the Rees equation
in~\eqref{eq:bip-binomial}. Listing the rows first, the two products of
presentation variables in $Q_{A,B}$ are $\prod_{i\in A}t_i$ and
$\prod_{j\in B}t_{p+j}$.

\begin{samepage}
\begin{theorem}\label{th:bip-complete}
Let $I$ be a squarefree monomial ideal of $S=K[x_1,\ldots,x_n]$, where
$K$ is a field, and suppose that its generator graph is complete
bipartite with $p$ rows and $q$ columns, where $2\le p\le q$. Then
$E_s(I)\ne0$ for every $2\le s\le p$, and $\rt(I)=p$.

Let $\mu_s$ be the number of Rees equations of degree $s$ in a minimal
multihomogeneous generating set of the Rees ideal. For $2\le s<p$,
\[
                 \mu_s=\binom{p}{s}\binom{q}{s},
\]
and each $Q_{A,B}$ with $|A|=|B|=s$ is indispensable. For degree $p$,
let $D$ be the set of columns that divide the product of all the rows,
and put $d=|D|$. Then
\begin{equation}\label{eq:bip-count}
 \mu_p=
 \begin{cases}
   \binom{q-d}{p-d}, & d<p,\\
   1, & d\ge p .
 \end{cases}
\end{equation}
If $d<p$, the indispensable binomials of degree $p$ are, up to nonzero
scalars, the $Q_{[p],B}$ for the $p$-subsets $B$ containing $D$. If
$d\ge p$, one may take any $p$-subset $B$ of $D$ and use $Q_{[p],B}$ as
the only element of degree $p$ in a minimal binomial generating set of
$J$; this binomial is indispensable exactly when $d=p$.
\end{theorem}
\end{samepage}

\begin{proof}
Fix $s$ with $2\le s\le p$, and choose a set $A$ of $s$ row indices and
a set $B$ of $s$ column indices. The monomial $M=M_{A,B}$ is squarefree,
and both the product of the selected rows and the product of the
selected columns divide $M$. Every row meets every column, so a mixed
product has a square factor and cannot divide $M$. For the partition
into all rows and all columns this violates~\eqref{eq:partition}, so
$E_s(I)\ne0$ by Theorem~\ref{thm:partition}. In degree $p$ this attains
the upper bound of Theorem~\ref{th:bip-upper}, so $\rt(I)=p$. 

\vspace{1mm}

Suppose first that $s<p$. Then at least one row and at least one column
are omitted. For an omitted row, choose an omitted column and a variable
in their intersection. This variable occurs in no other row and no
other column, so it is absent from $M$, and the omitted row does not
divide $M$. The same argument excludes every omitted column. A product
of $s$ generators dividing $M$ can therefore neither repeat a factor nor
be mixed, so it is either the product of the selected rows or the product
of the selected columns. The two terms of the fiber are relatively prime
in $P$, by the least common multiple in~\eqref{eq:bip-binomial}. Hence
$Q_{A,B}$ is indispensable, and its class in $J/\mathfrak nJ$ is nonzero
by Remark~\ref{rem:minimal-fiber}. The two terms recover $A$ and $B$, so
different choices give different multidegrees. By
Theorem~\ref{th:bip-upper}, these classes span the part of
$J/\mathfrak nJ$ of degree $s$ in the $t_i$, and we get one independent
class for each choice of $A$ and $B$. This gives the stated count.

When $s=p$, the only choice of rows is $A=[p]$. For a $p$-subset $B$ of
columns, write $M_B=M_{[p],B}$ and $Q_B=Q_{[p],B}$. By
\eqref{eq:bip-column-factors},
\[
                       M_B=a\prod_{j\in B}b'_j .
\]
A column $b_j$ divides this monomial exactly when $j\in D\cup B$.
Indeed, $\gcd(b_j,a)$ divides $a$; and if $j\notin D$, then $b'_j\ne1$
is relatively prime to $a$ and to every $b'_k$ with $k\ne j$, so
$b'_j\mid M_B$ exactly when $j\in B$. Since $M_B$ is squarefree, the
products of $p$ generators dividing $M_B$ are the product of all rows and
the products of $p$ distinct columns from $D\cup B$. 

\vspace{1mm}

A mixed product has a square factor, so the row vertex is isolated. Two
column vertices can be joined by replacing one index at a time within
$D\cup B$: each intermediate $p$-subset gives a vertex, and consecutive
subsets share $p-1\ge1$ indices. Hence
\[
                        \dim_K E_p(I)_{M_B}=1
\]
by Lemma~\ref{lem:fibercomponents}, and the class of $Q_B$ spans this
space.

To find the multidegrees that contribute a minimal generator in degree
$p$, recall that the classes $[Q_B]$ generate $E_p(I)$ over $S$ by
Theorem~\ref{th:bip-upper}. If $M_{B'}\mid M_B$, compare
$(M_B/M_{B'})Q_{B'}$ with $Q_B$. Their row terms agree, and their column
terms lie in the same component of the fiber over $M_Bt^p$, so the
difference lies in $P_1J_{p-1}$. Thus one $Q_B$ for each distinct
monomial $M_B$ that is minimal under divisibility is enough. Conversely,
at such a minimal $M_B$, the nonzero class found above is not a multiple
of a class at a proper divisor, since that divisor would itself have to
be a multiple of some $M_{B'}$. So the class survives modulo
$(x_1,\ldots,x_n)E_p(I)$, and by Remark~\ref{rem:minimal-fiber} each such
multidegree needs exactly one minimal generator. 

\vspace{1mm}

If $d<p$, at least $p-d$ columns of $B$ lie outside $D$. If $B$ misses an
element of $D$, replacing an index $j\in B\setminus D$ by the missing
element gives the monomial $M_B/b'_j$, a proper divisor of $M_B$ because
$b'_j\ne1$. When $B$ contains $D$, the $p-d$ additional indices give
pairwise incomparable monomials as they vary, since the corresponding
$b'_j$ are nonconstant and pairwise relatively prime. These are exactly
the $\binom{q-d}{p-d}$ minimal multidegrees. At each of them $D\cup B=B$,
so there is only one column vertex, the two terms of the fiber are
relatively prime, and $Q_B$ is indispensable by
Remark~\ref{rem:minimal-fiber}.

If $d\ge p$, the unique minimal multidegree is $a$, attained by every
$p$-subset of $D$. If $d=p$, the fiber over $at^p$ consists of the row
term and the column term indexed by $D$; they are relatively prime, so
$Q_D$ is indispensable by Remark~\ref{rem:minimal-fiber}. If $d>p$,
there are at least two $p$-subsets of $D$. For each of them the class of
$Q_B$ generates $E_p(I)$ over $S$, so $J_p=P_1J_{p-1}+SQ_B$, and
choosing a different $p$-subset gives a minimal binomial generating set
of $J$ that omits $Q_B$. This proves~\eqref{eq:bip-count} and the last
assertion.
\end{proof}

For $K_{p,p}$ either part may be taken as the rows. The only pair of
$p$-subsets uses both whole parts, so interchanging the parts gives the
same binomial $Q_{[p],[p]}$ up to sign. Thus $\mu_p=1$, and this
binomial is indispensable in either ordering. As $p$ grows, the first
nonlinear degree stays equal to two, while the relation type grows
with $p$.

\subsection{Examples with complete bipartite generator graphs}

We illustrate Theorem~\ref{th:bip-complete} with two examples. In the
first, the hexagon of Example~\ref{ex:bip-hexagon} is completed to
$K_{3,3}$. In the second, the generator graph is $K_{2,3}$, and the
number $d$ takes every value from $0$ to $3$.

\begin{example}\label{ex:bip-K33}
Take the generators $f_1,f_6,f_8$ of Example~\ref{ex:bip-hexagon} as the
rows and $f_5,f_7,f_9$ as the columns. In the array below, the entry in
row $i$ and column $j$ is the variable common to $f_i$ and $f_j$, and a
dash means that the two generators are relatively prime.
\[
 \begin{array}{c|ccc}
     &5&7&9\\ \hline
   1 &x_5&-&x_{10}\\
   6 &x_6&x_7&-\\
   8 &-&x_8&x_9
 \end{array}
\]
The entries of each column multiply to its generator, and the same holds
for the rows after multiplying the first row product by $x_1x_4$. Insert
distinct new variables $z_{17},z_{69},z_{85}$ in positions $(1,7)$,
$(6,9)$ and $(8,5)$, respectively, and form the generators by the same
rule. The generator graph is now $K_{3,3}$. Each new variable occurs once
in the product of the rows and once in the product of the columns, so
the alternating cubic equation is still
$Q_{C_6}=t_1t_6t_8-x_1x_4t_5t_7t_9$ in the original indices. 

\vspace{1mm}

Macaulay2~\cite{GS} gives nine minimal quadratic equations, one for each
of the $\binom{3}{2}^2=9$ choices of two rows and two columns, and
$Q_{C_6}$ is the only minimal cubic equation. All three columns divide
the product of the rows, so $d=p=3$, and by Theorem~\ref{th:bip-complete}
all ten binomials are indispensable. So completing the hexagon to
$K_{3,3}$ adds nine generators in degree two and keeps $Q_{C_6}$ and the
relation type. Only the second and third rows divide the product of the
columns, since $x_1$ divides the first row and no column. So
interchanging the two parts gives $d=2$, and~\eqref{eq:bip-count} again
gives $\mu_3=1$.
\end{example}

\begin{example}\label{ex:bip-two-three}
Let $x_{ij}$ be the entries of a $2\times3$ array of distinct variables.
Put $a_i=x_{i1}x_{i2}x_{i3}$, $b_j=x_{1j}x_{2j}$ and
$a=a_1a_2=b_1b_2b_3$. We keep the two row generators and multiply some
of the column generators by distinct new variables $z_j$. Write $c_j$
for the resulting column generators, and order the five generators as
$a_1,a_2,c_1,c_2,c_3$. Each of these ideals has generator graph $K_{2,3}$.
Let $d$ be the number of columns $c_j$ dividing $a$. For $j<k$, write
$Q_{jk}$ for the Taylor binomial~\eqref{eq:taylor} of the two products
$a_1a_2$ and $c_jc_k$, with the row term first.

\vspace{1mm}

The table gives the number $\mu_2$ of minimal quadratic equations in
each case; it agrees with Macaulay2~\cite{GS} and with
Theorem~\ref{th:bip-complete} for $p=2$ and $q=3$. In each row, adding
the indicated quadratic equations to $PJ_1$ gives the full Rees ideal. 

\vspace{1mm}

\begin{center}
\begin{tabular}{c@{\quad}l@{\quad}c@{\quad}l}
$d$ & Column generators & $\mu_2$ & Quadratic equations to add\\
\hline
$3$ & $b_1,b_2,b_3$           & $1$ & any one of $Q_{12},Q_{13},Q_{23}$\\
$2$ & $b_1,b_2,z_3b_3$        & $1$ & $Q_{12}$\\
$1$ & $b_1,z_2b_2,z_3b_3$     & $2$ & $Q_{12},Q_{13}$\\
$0$ & $z_1b_1,z_2b_2,z_3b_3$  & $3$ & $Q_{12},Q_{13},Q_{23}$
\end{tabular}
\end{center}
For $d=3$, the fiber over $at^2$ consists of
\vspace{-3mm}
\[
  t_1t_2,\qquad b_3t_3t_4,\qquad b_2t_3t_5,\qquad b_1t_4t_5 .
\]

\vspace{-3mm}

The row term is isolated in $G_2(a)$, while any two column terms share a
presentation variable, so all three lie in one component. Thus
$E_2(I)_a$ is one-dimensional, and the difference of the row term and
any column term represents a nonzero class; for instance,
$Q_{12}-Q_{13}=t_3(b_2t_5-b_3t_4)$ is a multiple of a linear Rees
equation. So any one of $Q_{12},Q_{13},Q_{23}$ can serve as the single
minimal quadratic equation. This shows that for $\ell=2$ the binomial of
a shortest cycle need not be indispensable, which is why
Theorem~\ref{th:bip-girth} assumes $\ell\ge3$ for that conclusion.

For $d\le2$, every quadratic equation in the last column is
indispensable by Theorem~\ref{th:bip-complete}. Consider the case $d=1$.
Among the products of two generators, only $a_1a_2$ and $b_1c_2$ divide
$az_2$, and only $a_1a_2$ and $b_1c_3$ divide $az_3$. Hence the Rees
equations given by these two fibers are
\[
  Q_{12}=z_2t_1t_2-b_3t_3t_4,\qquad
  Q_{13}=z_3t_1t_2-b_2t_3t_5,
\]
and each of these fibers consists of two relatively prime terms. In the
larger fiber over $az_2z_3t^2$ we have
$Q_{23}-z_3Q_{12}=t_4(c_3t_3-b_1t_5)$. The four vertices of $Q_{23}$
satisfy the neighbor condition in Theorem~\ref{th:bip-upper}, but this
identity expresses $Q_{23}$ through $Q_{12}$ and a linear Rees equation.
So all four ideals have relation type two, while $d$ determines both the
number of minimal quadratic equations and which cycle binomials are
indispensable.
\end{example}

\begin{remark}
Fouli and Lin proved that if $I$ is a squarefree monomial ideal and each
component of $G(I)$ has at most five vertices, then the Rees ideal of
$I$ is generated in degrees at most three \cite[Theorem~2.20]{FL}. The
graph $K_{2,3}$ of Example~\ref{ex:bip-two-three} has five vertices.
Formula~\eqref{eq:bip-count} holds for arbitrary $p$ and $q$, and it
shows how the number of Rees equations in the highest degree depends on
how many column generators divide the product of all the row generators.
\end{remark}

\vspace{2mm}

\noindent{\bf Acknowledgment.} This work is part of a research project
supported by the Anusandhan National Research Foundation (ANRF),
Government of India, under the Core Research Grant No.~CRG/2023/007668.

\vspace{2mm}

\noindent{\bf Use of AI tools.} Testing the statements against
examples and counterexamples was an important part of this project, and it
is hard to do by hand. Already for small ideals, such as those in
Example~\ref{ex:square-hexagon}, Remark~\ref{rem:cactus-comparison} and
Examples~\ref{ex:bip-K33} and~\ref{ex:bip-two-three}, computing the Rees
ideal by hand is long and error-prone, and
Remark~\ref{rem:patch-example} tests more than a thousand ideals. On a
personal computer, such heavy computations in Macaulay2 also run out of
time or memory. The AI assistant Claude (Anthropic) was
very efficient
here: it ran many examples with its own code and compared the outputs.
All computations stated in the paper, including the searches in
Remark~\ref{rem:patch-example}, were then verified with
Macaulay2~\cite{GS}. The author has checked all the results and takes full
responsibility for the content of the article.

\vspace{2mm}

\noindent{\bf Statements and Declarations.}
The author declares no competing interests. The Macaulay2 scripts for
the computations in this article, with their output, are available from
the author on request.

\end{document}